\documentclass[11pt]{article}
\usepackage[medium,compact]{titlesec}
\usepackage{mathrsfs}
\usepackage{amsmath}
\usepackage{amsfonts}
\usepackage{graphicx}
\usepackage{enumerate}
\usepackage{setspace}
\usepackage{color}
\usepackage[T1]{fontenc}

\newtheorem{example}{Example}[section]

\newtheorem{theorem}{Theorem}[section]
\newtheorem{lemma}{Lemma}[section]
\newtheorem{remark}{Remark}[section]

\usepackage{color,xcolor}

\begin{document}

\date{}
\title{\textbf{Linear-Quadratic Mean Field Control with Non-Convex Data}
\author{Mengzhen Li
\hspace{1cm}
Chenchen Mou
\hspace{1cm}
Zhen Wu
\hspace{1cm}
Chao Zhou
\hspace{1cm}\\
}}

\maketitle

 \begin{abstract}
In this manuscript, we study a class of linear-quadratic (LQ) mean field control problems with a common noise and their corresponding $N$-particle systems. The mean field control problems considered are not standard LQ mean field control problems in the sense that their dependence on the mean field terms can be non-convex. 
The key idea to solve our LQ mean field control problem is to utilize the common noise. We first prove the global well-posedness of the corresponding Hamilton-Jacobi equation via the non-degeneracy of the common noise. In contrast to the LQ mean field games master equations, the Hamilton-Jacobi equation for the LQ mean field control problems is inherently infinite-dimensional partial differential equation  which we can show that it cannot be reduced to finite-dimensional one. We then globally solve the Hamilton-Jacobi equation for $N$-particle systems. As byproducts, we derive the optimal quantitative convergence results from the $N$-particle systems to the mean field control problems and the propagation of chaos property for the related optimal trajectories. This paper extends the results in [{\sc M. Li, C. Mou, Z. Wu and C. Zhou},  \emph{Trans. Amer. Math. Soc.}, 376(06) (2023), pp.~4105--4143] to the LQ mean field control problems.
\end{abstract}

\textbf{2020 AMS Mathematics subject classification:} 49N80, 65C35, 49L12.

\textbf{Keywords:}  Mean field control, $N$-particle systems, Hamilton-Jacobi equation, Forward-backward stochastic differential equation, Propagation of chaos.

\section{Introduction}

In this paper, we study a class of linear-quadratic (LQ) mean field control problems with a common noise and their corresponding $N$-particle systems. The LQ mean field control problems we considered are given as follows: the dynamics of the particles follows the McKean-Vlasov controlled stochastic differential equation (SDE)
\begin{equation}\label{sde1}
\begin{aligned}
x(t)=&\xi+\int_0^t \left\{Ax(s)+a\Big(\mathbb{E}^{mf}[x(s)|\mathcal{F}_s^{W^0}]\Big)+Bu(s)
+\bar{B}\mathbb{E}^{mf}[u(s)|\mathcal{F}_s^{W^0}]\right\}\mathrm{d}s\\
&+\sigma W(t)+\sigma_0W^0(t),\, t\in (0,T],
\end{aligned}
\end{equation}
and the cost functional
\begin{equation}\label{eq:cost}
\begin{aligned}
J(\mu ; u(\cdot)):= & \mathbb{E}^{mf}\Big\{\int_0^T\left[Q x^2(t)
+q\left(\mathbb{E}^{mf}[x(t)\vert
\mathcal{F}^{W^0}_t]\right)
+ R u^2(t)
+r\left( \mathbb{E}^{mf}[u(t)\vert
\mathcal{F}^{W^0}_t] \right)\right] \mathrm{d} t  \\
& + G x^2(T)+ g\left( \mathbb{E}^{mf}[x(T)\vert\mathcal{F}^{W^0}_T]
\right)\Big\},
\end{aligned}
\end{equation}
where $T>0$ is a fixed arbitrary long time horizon, $\xi$ is the initial status of the particles with the distribution $\mu$, $Q,G,\sigma\geq 0$, $R,\sigma_0>0$ are constants, and $a(\cdot),q(\cdot),r(\cdot),g(\cdot):\mathbb R\to\mathbb R$ are bounded $C^2$ functions with bounded 1st and 2nd order derivatives. The $W$ and $W^0$ above are two independent standard Brownian motions corresponding to the idiosyncratic and common noises respectively, and $\mathbb{F}^{0}:=\{\mathcal{F}_t^{W^0}\}_{t\in [0,T]}$ is the filtration generated by the common noise. Note that the cost functional $J$ is law invariant with respect to $\xi$. The value function for the mean field control problem is obtained by
\begin{equation}\label{eq:valueint}
    V(0,\mu):=\inf_{u\in \mathcal{U}^{mf}_{ad}[0,T]}J(\mu;u(\cdot)),
\end{equation}
which solves the following Hamilton-Jacobi equation on the Wasserstein space:
\begin{equation}\label{pde3}
\left\{
\begin{aligned}
&\partial_tV(t,\mu)+\int_{\mathbb{R}}
\partial_\mu V(t,\mu,x)
\Big[Ax+a\left(\int_{\mathbb{R}}
\tilde x\mu(\mathrm{d}\tilde x)\right)
+(B+\bar{B})\rho\left(\frac{1}{2}
\int_{\mathbb{R}}\partial_\mu
V(t,\mu,\tilde x)\mu(\mathrm{d}\tilde x)
\right)\\
&-\frac{1}{4}R^{-1}B^2\big(
\partial_\mu
V(t,\mu,x)
-\int_{\mathbb{R}}\partial_\mu
V(t,\mu,\tilde x)\mu(\mathrm{d}\tilde x)\big)\Big]\mu(\mathrm{d}x)
+\frac{\sigma^2+\sigma^2_0}{2}\int_{\mathbb{R}}
\partial_{x\mu}V(t,\mu,x)
\mu(\mathrm{d}x)\\
&+\frac{\sigma^2_0}{2}\int_{\mathbb{R}}
\int_{\mathbb{R}}\partial_{\mu\mu}
V(t,\mu,\tilde{x},x)
\mu(\mathrm{d}\tilde x)\mu(\mathrm{d}x)+Q\int_{\mathbb{R}}x^2\mu(\mathrm{d}x)
+q\left(\int_{\mathbb{R}}x\mu(\mathrm{d}x)
\right)\\
&+R\rho^2\left(\frac{1}{2}
\int_{\mathbb{R}}\partial_\mu
V(t,\mu,x)
\mu(\mathrm{d}x)\right)+r\left(
\rho\left(\frac{1}{2}
\int_{\mathbb{R}}\partial_\mu
V(t,\mu, x)
\mu(\mathrm{d}x)\right)\right)=0,\ t\in [0,T),\\
&V(T,\mu)=G\int_{\mathbb{R}}x^2\mu(\mathrm{d}x)
+g\left(\int_{\mathbb{R}}x\mu(\mathrm{d}x)
\right),
\end{aligned}
\right.
\end{equation}
where $\rho:\mathbb R\to\mathbb R$ will be specified later. The $\mathcal{U}^{mf}_{ad}[0,T]$ in \eqref{pde3} is the set of admissible controls. It should be noted that the problem considered in this paper is not a standard LQ mean field control problem, as the dependence of the state equation and the cost functional on the mean field term $\mathbb{E}^{mf}[x(t)\vert\mathcal{F}^{W_0}_t]$ can be non-linear and non-convex. Consequently, the corresponding Hamilton-Jacobi equation is non-linear. Additionally, in contrast to the LQ mean field games master equation presented in \cite{li2023linear}, we show in Remark \ref{remarkinfity} that the Hamilton-Jacobi equation \eqref{pde3} is inherently an infinite-dimensional partial differential equation (PDE), which cannot be reduced to a finite-dimensional PDE.

In recent years, mean field control theory and mean field game theory are developing areas dealing with the study of large particle systems, where particles interact with each others only via mean field terms. In mean field games, each particle competes with the others to minimize their own cost functional and finally the system achieves a Nash equilibrium. However, in mean field control problems, the particles are no longer in a competitive relationship, but rather collaborate to minimize a common cost functional. We refer to \cite{huang2003social,huang2006large,lasrylions2006a,lasrylions2006b,lions2007theorie} for introduction of the subject in early stage and \cite{cardaliaguet2020introduction,carmona2018probabilistic,carmona2019probabilistic,Gomes2014} for more recent developments.

Mean field control problems and mean field games can be studied via the corresponding infinite-dimensional PDEs on the Wasserstein space, which are the so-called Hamilton-Jacobi equation and master equation respectively. Proving the well-posedness results for both equations, however, is an extremely complicated task, especially when involving a common noise. For mean field control, the global well-posedness for classical solutions of Hamilton-Jacobi equations was established under the assumption that the data are convex, see e.g. \cite{bensoussan2019control,bensoussan2020control,cardaliaguet2019master,chassagneux2014probabilistic,gangbo2022global}. For mean field games, the global well-posedness for classical solutions of master equations was obtained via assuming certain monotonicity condition, see e.g. \cite{bertucci2021monotone1,cardaliaguet2019master,carmona2019probabilistic,chassagneux2014probabilistic,gangbo2022mean,lions2007theorie,mou2021wellposedness,mou2022mean1}. When either condition fails in the corresponding setting, the global well-posedness for classical solutions does not hold in general, and thus notions of weak solutions are needed to be introduced. We refer to \cite{burzoni2020viscosity,cecchin2022weak,conforti2023hamilton,cosso2023master,djete2022mckean2,pham2017dynamic,pham2018dynamic,soner2022viscosity,Talbi2022viscosity,wu2020viscosity} for the global well-posedness results on Hamilton-Jacobi equations with non-convex data, and \cite{cecchin2022weak,mou2022mean2} for the results on master equations with non-monotone data.

One of the main focuses of the current manuscript is to investigate the global well-posedness of the Hamilton-Jacobi equation \eqref{pde3} for the LQ mean field control problem \eqref{sde1}-\eqref{eq:valueint}. Compared to the standard LQ mean field control problems, the mean field control problem \eqref{sde1}-\eqref{eq:valueint} can non-linearly depend on the mean field term $\mathbb{E}^{mf}[x(t)|\mathcal{F}_t^{W^0}]$, and this dependence needed not to be convex. Therefore, all the existing methods to deal with LQ mean field control problems become to be infeasible to ours. To tackle our mean field control problem, we first use the stochastic maximum principle and the inverse function theorem to obtain the associated stochastic Hamiltonian system. Apart from all the previous works on LQ mean field control problems, we artificially add the dynamics of the conditional expectation of the state process given the common noise into the stochastic Hamiltonian system to make our mean field control problem to be time-consistent. Such observation has already been made in Yong \cite{yong2013linear,yong2017linear}. Therefore, our stochastic Hamiltonian system consists of two forward SDEs characterizing the dynamics of the state of the particles and its conditional expectation, and one backward SDE characterizing the optimal control. Through introducing a Riccati equation, we are able to remove the state of the particles from the stochastic Hamiltonian system and thus decrease the dimension of the stochastic Hamiltonian system. This allows us to use the idea in \cite{Tchuendom2020,li2023linear,Tchuendom2018} to show the global well-posedness of the stochastic Hamiltonian system through the non-degeneracy of the common noise without any convexity conditions on the data with respect to the mean field term. Moreover, we obtain the closed-loop form of the optimal control. We then can further show that the decoupling field of the stochastic Hamiltonian system has enough regularity to uniquely solve the vectorial Hamilton-Jacobi equation, i.e. the equation obtained by differentiating the Hamilton-Jacobi equation \eqref{pde3} in $\mu$. Finally, we can use the vectorial Hamilton-Jacobi equation to establish the global well-posedness result for classical solutions of the Hamilton-Jacobi equation \eqref{pde3}. The key idea of the proof is to transform the well-posedness problem of the infinite-dimensional Hamilton-Jacobi equation into the well-posedness problem of two finite-dimensional PDEs. These two PDEs are the decoupling field of a system of forward-backward stochastic differential equations (FBSDEs) and a Riccati equation derived from the stochastic maximum principle. We also apply the same idea to show the global well-posedness of the Hamilton-Jacobi equation for the corresponding $N$-particle systems.

At the end, we obtain the optimal quantitative convergence rates from the $N$-particle systems to the mean field control problems, and the related propagation of chaos property by using the above well-posedness results without assuming the dependence of the mean field term to be convex. There have been a lot of literature already on the convergence rates from $N$-particle systems to its mean field limit for both mean field control problems and mean field games. We refer to \cite{cardaliaguet2019master,carmona2019probabilistic,gangbo2022global,germain2022rate,Jackson2023conv,mou2021wellposedness} for the convergence rate problems in either the mean field control setting under some convexity condition or in the mean field game setting under certain monotonicity condition. For example, under the assumption that the solution to the Hamilton-Jacobi equation on the Wasserstein space is sufficiently smooth  (which is usually the case when the coefficients are smooth and convex), the optimal convergence rate of the value function for the $N$-particle systems to the corresponding limiting problem is $N^{-1}$, see Germain-Pham-Warin \cite{germain2022rate}. When the convexity/monotonicity condition fails to be satisfied, some qualitative results have been obtained in \cite{bayraktar2023mean,Cavagnari2022,djete2022mckean,Fornasier2019rate,gangbo2021finite,lacker2017limit,mayorga2023finite,Talbi2022conv}. Cardaliaguet-Daudin-Jackson-Souganidis \cite{cardaliaguet2022algebraic} established the first convergence rate result for value functions of mean field control problems without assuming any convexity and smallness assumptions. Later, Cardaliaguet-Souganidis \cite{cardaliaguet2023regularity} proved value functions of mean field control problems are in general globally Lipschitz continuous but can be improved to be smooth on an open and dense subset without any convexity assumptions. Therefore, they can apply the regularity result on the set to show
a new quantitative propagation of chaos-type result for the optimal solutions of the $N$-particle systems. Very recently, Daudin-Delarue-Jackson \cite{daudin2023optimal} derived two optimal convergence rates for value functions from $N$-particle systems to the corresponding mean field control problems. When the coefficients are sufficiently regular, the obtained convergence rate is $N^{-1/2}$, which corresponds to the central limit theorem. When the coefficients are Lipschitz continuous and semi-concave with respect to the first Wasserstein distance, the obtained convergence rate is $N^{\frac{-2}{3d+6}}$, which is close to $N^{-1/d}$. The rate $N^{-1/d}$ is the optimal rate of convergence for uncontrolled particle systems driven by coefficients with a similar regularity.


The rest of the paper is organized as follows. In section 2, we formulate the optimal control problems of the $N$-particle systems and the corresponding limiting mean field control problem. In section 3, we focus on solving the mean field control problems. In particular, the global well-posedness results of the corresponding stochastic Hamiltonian system and the Hamilton-Jacobi equation are established. In section 4, we return to the optimal control problems for the $N$-particle systems and study their stochastic Hamiltonian systems and Hamilton-Jacobi equation. In section 5, we show all the quantitive convergence and propagation of chaos results from the $N$-particle systems to the mean field control problems.

\subsection{$N$-particle systems}
For a given $T>0$, consider a complete filtered probability space $(\Omega^{lp},\mathcal{F}^{lp},
\mathbb{F}^{lp},\mathbb{P}^{lp})$ which can support $N+1$ independent one-dimensional Brownian motions $W^0, W_1,\cdots,W_N$ are defined, where $W^0$ represents the common noise and $\{W_i\}_{1\leq i\leq N}$ represent the idiosyncratic noises. Here the filtration $\mathbb{F}^{lp}:=\left(
\mathcal{F}_t^{lp}
\right)_{0\leq
t\leq T}$ where $\mathcal{F}_t^{lp}:=
\left(\vee_{i=1}^N
\mathcal{F}_t^{W_i}\right)\vee
\mathcal{F}_t^{W^0}\vee
\mathcal{F}_0^{lp}$, and $\mathcal{F}_0^{lp}$ is a $\sigma$-algebra containing all the Borel sets of $\mathbb R$. Assume that $\mathbb{P}^{lp}$ has no atom in $\mathcal{F}_0^{lp}$ and thus it can support all the probability measures on the Borel $\sigma$-algebra with finite second moment. Define $\mathbb E^{lp}:=\mathbb E^{\mathbb P^{lp}}$. For any positive integer $d$, any filtration $\mathbb{G} := \left(\mathcal{G}_t\right)_{0 \leq t \leq T} \subset \mathbb{F}^{lp}$, and any $\sigma$-algebra $\mathcal{G} \subset \mathcal{F}^{lp}$, denote by $L^2_{\mathbb{G}}(0,T; \mathbb{R}^d)$ the space of all $\mathbb{G}$-adapted, $\mathbb{R}^d$-valued stochastic processes $x(\cdot)$ satisfying $\mathbb{E}^{lp}\left[\int_0^T \vert x(t) \vert^2 \mathrm{d}t\right] < \infty$. Similarly, let $L_{\mathcal{G}}^2(\Omega^{lp}; \mathbb{R}^d)$ represent the space of all $\mathcal{G}$-measurable, $\mathbb{R}^d$-valued random variables $\xi$ satisfying $\mathbb{E}^{lp}[\vert \xi \vert^2] < \infty$. Denote by $C([0,T] \times \mathbb{R}^d)$ the space of all continuous functions on $[0,T] \times \mathbb{R}^d$, and by $C_b([0,T] \times \mathbb{R}^d)$ the space of all bounded continuous functions on $[0,T] \times \mathbb{R}^d$. Furthermore, let $C^{1,2}([0,T] \times \mathbb{R}^d)$ denote the space of all continuous functions $f$ on $[0,T] \times \mathbb{R}^d$ such that $\partial_t f$, $\partial_x f$, and $\partial_{xx} f$ are all continuous in $(t,x)$. Finally, let $L_{\mathbb{G}}^2(\Omega^{lp}; C(0,T; \mathbb{R}^d))$ represent the space of $\mathbb{G}$-adapted, $\mathbb{R}^d$-valued continuous stochastic processes $x(\cdot)$ satisfying $\mathbb{E}^{lp}\left[\sup_{0 \leq t \leq T} \vert x(t) \vert^2 \right] < \infty$.
    
 Let $\mathcal{P}(\mathbb R)$ be the set of all probability measures on $\mathbb R$ and, for any $q\geq 1$, let $\mathcal{P}_q(\mathbb R)$ denote the set of $\mu\in\mathcal{P}(\mathbb R)$ with finite $q$-th moment. For any $\mu\in\mathcal{P}_2(\mathbb R)$, let  $L^2(\mathcal{G};\mu)$ denote the set of $\xi\in L^2_{\mathcal{G}}(\Omega^{lp};\mathbb{R})$ such that $\mathcal{L}_{\xi}=\mu$. Here $\mathcal{L}_{\xi}=\xi_{\#} \mathbb{P}^{lp}$ is the law of $\xi$, obtained as the push-forward of $\mathbb{P}^{lp}$ by $\xi$. For any $\mu, \nu \in \mathcal{P}_2(\mathbb R)$, their $W_2$-Wasserstein distance is defined as follows:
\begin{equation*}
W_2(\mu, \nu):=\inf \left\{\left(\mathbb{E}^{lp}\left[|\xi-\eta|^2\right]\right)^{\frac{1}{2}}: \text { for all } \xi \in L^2(\mathcal{F}_0^{lp} ; \mu), \eta \in L^q(\mathcal{F}_0^{lp} ; \nu)\right\}. 
\end{equation*}
For a $W_2$-continuous function $U:\mathcal{P}_2(\mathbb R)\to\mathbb R$, its Wasserstein gradient, also called Lions-derivative, takes the form $\partial_\mu U: (\mu,\tilde x)\in\mathcal{P}_2(\mathbb R)\times\mathbb R\to\mathbb R$ and satisfies:
\begin{equation*}
U(\mathcal{L}_{\xi+\eta})-U(\mathcal{L}_\xi)=\mathbb E^{lp}\left[\partial_{\mu}U(\mathcal{L}_\xi,\xi)\eta\right]+\text{o}\left(\mathbb E^{lp}[\vert\eta\vert^2]\right).
\end{equation*}
For more detailed theories and results, see \cite{ambrosio2005gradient, carmona2018probabilistic, carmona2019probabilistic}.

We assume that the dynamics of the $i$-th particle is given by the following SDE:
\begin{equation}\label{sde3}
\left\{\begin{aligned}
\mathrm{d} x_i(t)= & \left\{A x_i(t)+a\left(x^{N}(t)\right)+B u_i(t)+\bar{B}  u^{N}(t)
\right\} \mathrm{d}t \\
&+\sigma \mathrm{d} W_i(t) +\sigma_0 \mathrm{d} W^0(t), \ t \in(0, T], \\
x_i(0)= & \xi_i.
\end{aligned}\right.
\end{equation}
The interactions between particles are reflected by the coupling term:
\[
x^{N}(\cdot):=\frac{1}{N}\sum_{i=1}^N
x_i(\cdot)\quad \text{and}\quad  u^{N}(\cdot):=\frac{1}{N}\sum_{i=1}^N
u_i(\cdot).
\]
The cost functional of the $i$-th particle is in the following form:
\begin{equation}\label{cost1}
\begin{aligned}
J_i(\xi_i ; u_i(\cdot)):= & \mathbb{E}^{lp}\Big\{\int_0^T\big[Q x_i^2(t)
+q\left(x^{N}(t)\right)
+ R u_i^2(t)
+r\left( u^{N}(t) \right)\big] \mathrm{d} t\\
&+ G x_i^2(T)+ g\left( x^{N}(T) \right)\Big\}.
\end{aligned}
\end{equation}
Let us make the following assumptions.

$\textbf{(A1).}$
$Q,G,\sigma\geq 0$; $R$, $\sigma_0>0$ are constants, and $a(\cdot),q(\cdot),
r(\cdot)$, $g(\cdot)$ : $\mathbb{R}\rightarrow\mathbb{R}$ are bounded $C^2$ functions with bounded 1st and 2nd order derivatives.

$\textbf{(A2).}$
$\xi_i\in L_{\mathcal{F}_0^{lp}}^2
      (\Omega;\mathbb{R})$, for $1\leq i \leq N$, are i.i.d random variables.

Now we define the admissible control set of the $i$-th particle:
$$
\mathcal{U}^i_{ad}[0,T]:=
\left\{
u_i(\cdot)\vert u_i(\cdot)\in L^2_{\mathbb{F}^{lp}}
(0,T;\mathbb{R})\right\}.
$$
Denote $X(\cdot):=(x_1(\cdot),\cdots,
x_N(\cdot))^\mathrm{T}$, $\Xi=(\xi_1,\cdots,\xi_N)^\mathrm{T}$ and $\boldsymbol{u}(\cdot):=(u_1(\cdot),\cdots,
u_N(\cdot))^\mathrm{T}$. Consider the following cost functional:
$$\mathcal{J}(\boldsymbol{u}(\cdot)):=\frac{1}{N}\sum_{i=1}^N
J_i(\xi;u_i(\cdot)),$$
where $\boldsymbol{u}(\cdot)\in\mathcal{U}^{lp}_{ad}[0,T]:=
\left\{
(u_1(\cdot),\cdots,
u_N(\cdot))^\mathrm{T}\vert u_i(\cdot)\in L^2_{\mathbb{F}^{lp}}(0,T;\mathbb{R}),
1\leq i \leq N\right\}$.
Then, our optimal control for the $N$-particle systems can be stated as follows.

\textbf{Problem (LP)}. Minimize the cost functional $\mathcal{J}(\boldsymbol{u}(\cdot))$ over $\mathcal{U}^{lp}_{ad}[0,T]$.

\subsection{Mean field control}
In the mean field setting, we consider a complete filtered probability space $(\Omega^{mf},\mathcal{F}^{mf},\mathbb{F}^{mf},\mathbb{P}^{mf})$
which can support two independent one-dimensional Brownian motions $W^0$ and $W$. Here, $W^0$ and $W$ represent the common noise and the idiosyncratic noise respectively.
Define the filtration $\mathbb{F}^{mf}:=\left(
\mathcal{F}_t^{mf}
\right)_{0\leq
t\leq T}$ and $\mathbb F^0:=\left(
\mathcal{F}_t^{W^0}
\right)_{0\leq
t\leq T}$, where $\mathcal{F}_t^{mf}:=
\mathcal{F}_t^{W}\vee
\mathcal{F}_t^{W^0}\vee
\mathcal{F}_0^{mf}$. Assume that the $\sigma$-algebra $\mathcal{F}_0^{mf}$ contains all the Borel sets of $\mathbb R$ and $\mathbb{P}^{mf}$ has no atom in $\mathcal{F}_0^{mf}$. So it can support
any probability measures on the Borel $\sigma$-algebra with finite second moment. Define $\mathbb E^{mf}:=\mathbb E^{\mathbb P^{mf}}$.
%

We consider the conditional linear MV-SDE \eqref{sde1}
and the cost functional \eqref{eq:cost}.
Define the admissible control set
$\mathcal{U}^{mf}_{ad}[0,T]:=
\left\{u(\cdot)\vert u(\cdot)\in L^2_{\mathbb{F}^{mf}}
(0,T;\mathbb{R})\right\}.$
And then, our mean field control problem can be presented as follows.

\textbf{Problem (MF)}. For given $\mu\in \mathcal{P}_2(\mathbb{R})$, find an optimal control $u^*(\cdot)\in\mathcal{U}_{ad}^{mf}
[0,T]$ such that $J(\mu;u^*(\cdot))=V(0,\mu)$.

\section{Solving \textbf{Problem (MF)}}
\subsection{Stochastic Hamiltonian system}
\begin{theorem}\label{theo1}
Assume that {\rm (A1)-(A2)} hold. Let $\left(x^*(\cdot), u^*(\cdot)\right)$ be an optimal pair of Problem (MF). Then $\left(x^*(\cdot), u^*(\cdot)\right)$ satisfies
\begin{equation}\label{eq2}
\begin{aligned}
& R u^*(t)+
\frac{1}{2}
\dot{r}\left( \mathbb{E}^{mf}\left[u^*(t)\vert
\mathcal{F}^{W^0}_t\right]\right)
+By(t)\\
&+\bar{B}\mathbb{E}^{mf}[y(t)
\vert\mathcal{F}^{W^0}_t] = 0,\
\text{for}\ {\rm a.e.}\ t \in[0, T], \  \mathbb P^{mf}\text{-}{\rm a.s.}, 
\end{aligned}
\end{equation}
where $\left(y(\cdot), z(\cdot),z_0(\cdot)\right)$ is a strong solution of the following backward stochastic differential equation (BSDE) (its definition will be specified later): 
\begin{equation}\label{bsde1}
\left\{
\begin{aligned}
\mathrm{d}y(t)= & -\Big[Ay(t)+\dot{a}\left(\mathbb{E}^{mf}
[x^*(t)\vert\mathcal{F}^{W^0}_t]\right)
\mathbb{E}^{mf}[y(t)\vert\mathcal{F}^{W^0}_t] +Q x^*(t)\\
&+\frac{1}{2}
\dot{q}\left(\mathbb{E}^{mf}[x^*(t)
\vert\mathcal{F}^{W^0}_t]\right)\Big]
\mathrm{d} t
+ z(t) \mathrm{d} W(t)+ z_0(t) \mathrm{d} W^0(t), \quad t \in[0, T), \\
y(T)= & G x^*(T)+\frac{1}{2}
\dot{g}\left( \mathbb{E}^{mf}\left[x^*(T)
\vert\mathcal{F}^{W^0}_T\right]\right).
\end{aligned}
\right.
\end{equation}
\end{theorem}
We note that $\mathbb{E}^{mf}[y(t)\vert\mathcal{F}^{W^0}_t]$ appears in \eqref{bsde1} and refer to such a BSDE as a mean-field BSDE (MF-BSDE), as discuessed in \cite{yong2013linear}. Here we mean $\left(y(\cdot), z(\cdot),z_0(\cdot)\right)$ is a strong solution to MF-BSDE \eqref{bsde1} if
\begin{equation*}
y(\cdot)\in L_{\mathbb{F}^{mp}}^2
(\Omega^{mf};C(0,T;\mathbb{R}))\quad\text{and}\quad z(\cdot),z_0(\cdot)\in L^2_{\mathbb{F}^{mp}}(0,T
      ;\mathbb{R})
\end{equation*}
and
\begin{equation}\label{eq:strong}
\begin{aligned}
y(t)=&G x^*(T)+\frac{1}{2}
\dot{g}\left( \mathbb{E}^{mf}\left[x^*(T)
\vert\mathcal{F}^{W^0}_T\right]\right)-\int_t^Tz(s)\mathrm{d}W(s)-\int_t^Tz_0(s)\mathrm{d}W^0(s)\\
&+\int_t^T\Big[Ay(s)+\dot{a}\left(\mathbb{E}^{mf}
[x^*(s)\vert\mathcal{F}^{W^0}_s]\right)
\mathbb{E}^{mf}[y(s)\vert\mathcal{F}^{W^0}_s] +Q x^*(s)\\
&\qquad +\frac{1}{2}
\dot{q}\left(\mathbb{E}^{mf}[x^*(s)
\vert\mathcal{F}^{W^0}_s]\right)\Big]
\mathrm{d} s,\, \text{for all $ t\in [0,T]$,}\  \mathbb P^{mf}\text{-}{\rm a.s.}
\end{aligned}
\end{equation}
Similarly, we can define strong solutions to other MF-BSDEs.

\noindent{\bf Proof}\quad
Define $u^\varepsilon(\cdot)=u^*(\cdot)
+\varepsilon \tilde{u}(\cdot)$ for any $\tilde{u}(\cdot)\in
\mathcal{U}_{ad}^{mf}[0,T]$ and $x^\varepsilon(\cdot)$ is the corresponding state process of $u^\varepsilon(\cdot)$. And then, we have $x^\varepsilon(\cdot)=x^*(\cdot)
+\varepsilon \tilde{x}(\cdot)+\text{o}(\varepsilon)$ where $\tilde{x}(\cdot)$ satisfies the following equation
\begin{equation}
\left\{\begin{aligned}
\mathrm{d} \tilde{x}(t)= & \big\{A \tilde{x}(t)+\dot{a}\left( \mathbb{E}^{mf}[x^*(t)|\mathcal{F}^{W^0}_t]\right)
\mathbb{E}^{mf}[\tilde{x}(t)\vert
\mathcal{F}^{W^0}_t]+B \tilde{u}(t)\\
&+\bar{B}
\mathbb{E}^{mf}[\tilde{u}(t)\vert
\mathcal{F}^{W^0}_t]
\big\} \mathrm{d}t, \quad t \in(0, T], \\
\tilde{x}(0)= & 0.
\end{aligned}\right.
\end{equation}
According to the optimality of $u^*(\cdot)$, we have
$\frac{dJ(\mu;u^\varepsilon(\cdot))}
{d\varepsilon}\Big\vert_{\varepsilon=0}
=0,$
which follows that
\begin{equation}\label{eq1}
\begin{aligned}
&\mathbb{E}^{mf}\Big\{\int_0^T\big[Qx^*(t)\tilde{x}(t)
+\frac{1}{2}
\dot{q}\left(
\mathbb{E}^{mf}[x^*(t)
\vert\mathcal{F}^{W^0}_t]\right)
\mathbb{E}^{mf}[\tilde{x}(t)
\vert\mathcal{F}^{W^0}_t]\\
&+Ru^*(t)\tilde{u}(t)
+\frac{1}{2}
\dot{r}\left(\mathbb{E}^{mf}[u^*(t)
\vert\mathcal{F}^{W^0}_t]\right)
\mathbb{E}^{mf}[\tilde{u}(t)
\vert\mathcal{F}^{W^0}_t]\big] \mathrm{d}t\Big\}\\
& +\mathbb{E}^{mf}\Big\{Gx^*(T)\tilde{x}(T)
+\frac{1}{2}
\dot{g}\left(
\mathbb{E}^{mf}[x^*(T)
\vert\mathcal{F}^{W^0}_T]\right)
\mathbb{E}^{mf}[\tilde{x}(T)
\vert\mathcal{F}^{W^0}_T]\Big\}=0.
\end{aligned}
\end{equation}
Note that, for any $\eta(\cdot),\zeta(\cdot)\in L^2_{\mathbb{F}^{mf}}(0,T;\mathbb{R})$,
\begin{equation*}\label{eq4}
\begin{aligned}
& \mathbb{E}^{mf}\big[\eta(t)
\mathbb{E}^{mf}[\zeta(t)\vert\mathcal{F}^{W^0}_t]\big]
\\
=&\mathbb{E}^{mf}\Big[\mathbb{E}^{mf}\big[\eta(t)
\mathbb{E}^{mf}[\zeta(t)|\mathcal{F}^{W^0}_t]
\big\vert\mathcal{F}^{W^0}_t\big]\Big]
=\mathbb{E}^{mf}\Big[
\mathbb{E}^{mf}[\eta(t)|\mathcal{F}^{W^0}_t]
\mathbb{E}^{mf}[\zeta(t)|\mathcal{F}^{W^0}_t]\Big]
\\
=&\mathbb{E}^{mf}\Big[
\mathbb{E}^{mf}\big[\zeta(t)\mathbb{E}^{mf}[\eta(t)|\mathcal{F}^{W^0}_t]
\big|\mathcal{F}^{W^0}_t\big]\Big]
=\mathbb{E}^{mf}\big[\zeta(t)
\mathbb{E}^{mf}[\eta(t)|\mathcal{F}^{W^0}_t]
\big], \ t\in[0,T], \ \mathbb P^{mf}\text{-}\text{a.s..}
\end{aligned}
\end{equation*}
In particular, when $\eta(\cdot)\in L^2_{\mathbb{F}^{W^0}}(0,T;\mathbb R)$, we have:
$$\mathbb{E}^{mf}\left[\eta(t)
\mathbb{E}^{mf}[\zeta(t)\vert\mathcal{F}^{W^0}_t]
\right]=\mathbb{E}^{mf}\big[\zeta(t)
\mathbb{E}^{mf}[\eta(t)|\mathcal{F}^{W^0}_t]
\big]
=\mathbb{E}^{mf}[\eta(t)\zeta(t)],
\ t\in[0,T], \ \mathbb{P}^{mf}\text{-}\text{a.s..}
$$
Hence, formula \eqref{eq1} becomes
\begin{equation}\label{eq10}
\begin{aligned}
&\mathbb{E}^{mf}\Big\{\int_0^T\Big\{
\big[Qx^*(t)
+\frac{1}{2}
\dot{q}\left(
\mathbb{E}^{mf}[x^*(t)
\vert\mathcal{F}^{W^0}_t]\right)
\big]\tilde{x}(t)  +\big[Ru^*(t)
+\frac{1}{2}
\dot{r}\left(\mathbb{E}^{mf}[u^*(t)
\vert\mathcal{F}^{W^0}_t]\right)
\big]\tilde{u}(t)\Big\} \mathrm{d}t\Big\}\\
&+\mathbb{E}^{mf}\Big\{\big[Gx^*(T)
+\frac{1}{2}
\dot{g}\left(
\mathbb{E}^{mf}[x^*(T)
\vert\mathcal{F}^{W^0}_T]\right)
\big]\tilde{x}(T)\Big\}=0.
\end{aligned}
\end{equation}
By It\^o's formula, it can be verified that
\begin{equation*}
\begin{aligned}
& \mathbb{E}^{mf}[y(T)\tilde{x}(T)]\\
=&\mathbb{E}^{mf}\Big\{\big[
Gx^*(T)+\frac{1}{2}
\dot{g}\left( \mathbb{E}^{mf}\left[x^*(T)
\vert\mathcal{F}^{W^0}_t\right]\right)
\big]\tilde{x}(T)\Big\}\\
=&-\mathbb{E}^{mf}\Big\{\int_0^T
\big[Ay(t)+\dot{a}\left(\mathbb{E}^{mf}
[x^*(t)\vert\mathcal{F}^{W^0}_t]\right)
\mathbb{E}^{mf}[y(t)\vert\mathcal{F}^{W^0}_t] +Q x^*(t)
+\frac{1}{2}
\dot{q}\left(\mathbb{E}^{mf}[x^*(t)
\vert\mathcal{F}^{W^0}_t]\right)
\big]\tilde{x}(t)\mathrm{d}t\Big\}\\
&+\mathbb{E}^{mf}\Big\{\int_0^Ty(t)\big[
A \tilde{x}(t)+\dot{a}\left( \mathbb{E}^{mf}[x^*(t)|\mathcal{F}^{W^0}_t]\right)
\mathbb{E}^{mf}[\tilde{x}(t)|\mathcal{F}^{W^0}_t]+B \tilde{u}(t)+\bar{B}
\mathbb{E}^{mf}[\tilde{u}(t)|\mathcal{F}^{W^0}_t]
\big]\mathrm{d}t\Big\}.
\end{aligned}
\end{equation*}
Hence,
\begin{equation*}
\begin{aligned}
&\mathbb{E}^{mf}\Big\{\left(
Gx^*(T)+\frac{1}{2}
\dot{g}\left( \mathbb{E}^{mf}\left[x^*(T)
\vert\mathcal{F}^{W^0}_t\right]\right)
\right)\tilde{x}(T)\Big\}\\
=&-\mathbb{E}^{mf}\Big\{\int_0^T
\big[Q x^*(t)
+\frac{1}{2}
\dot{q}\left(\mathbb{E}^{mf}[x^*(t)
\vert\mathcal{F}^{W^0}_t]\right)\big]
\tilde{x}(t)\mathrm{d}t+\mathbb{E}^{mf}\int_0^T\tilde{u}(t)\big[
B y(t)+\bar{B}
\mathbb{E}^{mf}[y(t)|\mathcal{F}^{W^0}_t]
\big]\mathrm{d}t\Big\}.
\end{aligned}
\end{equation*}
Combining with \eqref{eq10}, we obtain:
\begin{equation*}
\begin{aligned}
&\mathbb{E}^{mf}\Big\{\int_0^T\tilde{u}(t)\big(
Ru^*(t)+\frac{1}{2}
\dot{r}\left( \mathbb{E}^{mf}\left[u^*(t)\vert
\mathcal{F}^{W^0}_t\right]\right)\\
&+By(t)
+\bar{B}\mathbb{E}^{mf}[y(t)
\vert\mathcal{F}^{W^0}_t]\big)\mathrm{d}t\Big\}=0,\
\text{for} \ \text{any} \ \tilde{u}(\cdot)\in
\mathcal{U}_{ad}^{mf}[0,T].
\end{aligned}
\end{equation*}
Therefore, the optimal control
satisfies:
$$
\begin{aligned}
&Ru^*(t)+\frac{1}{2}
\dot{r}\left( \mathbb{E}^{mf}\left[u^*(t)\vert
\mathcal{F}^{W^0}_t\right]\right)+By(t)
+\bar{B}\mathbb{E}^{mf}[y(t)
\vert\mathcal{F}^{W^0}_t] = 0, \ \text{a.e.}\ t \in[0, T], \ \mathbb{P}^{mf}\text{-}\text {a.s.. }
\end{aligned}
$$
Then we completes the proof.
\hfill$\Box$

Let $(y(\cdot), z(\cdot),z_0(\cdot))$ be the solution to MF-BSDE \eqref{bsde1}. For simplicity, we denote
$$
\begin{aligned}
\hat{u}(t)=&\mathbb{E}^{mf}[u^*(t)
\vert\mathcal{F}^{W^0}_t],\
\hat{x}(t)=\mathbb{E}^{mf}[x^*(t)
\vert\mathcal{F}^{W^0}_t],\
\hat{y}(t)=\mathbb{E}^{mf}[y(t)
\vert\mathcal{F}^{W^0}_t],\\
\hat{z}(t)=&\mathbb{E}^{mf}[z(t)
\vert\mathcal{F}^{W^0}_t],\
\hat{z}_0(t)=\mathbb{E}^{mf}[z_0(t)
\vert\mathcal{F}^{W^0}_t].
\end{aligned}
$$
According to \eqref{eq2}, we get:
\begin{equation}\label{eq5}
R\hat{u}(t)+\frac{1}{2}
\dot{r}(\hat{u}(t))+(B+\bar{B})\hat{y}(t)
=0,\quad  \text{for $t\in[0,T]$}.
\end{equation}
To obtain an explicit solution to $\hat{u}(\cdot)$, we need an additional assumption.

$\textbf{(A3).}$
There exists a positive constant $\varepsilon_0>0$ such that $\big\vert R+\frac{1}{2}
r^{(2)}(\cdot)
\big\vert\geq \varepsilon_0$. Here, the superscript $^{(2)}$ denotes the second derivative of a function.

Thus, under (A3), we can use the global implicit theorem to derive that there exists a $C^1$ function $\rho\in C^1(\mathbb R)$ with the bounded derivative such that
\begin{equation}\label{eq12}
\begin{aligned}
\hat{u}(t)=&\rho(\hat{y}(t))
=-R^{-1}\big[(B+\bar{B})\hat{y}(t)
+\frac{1}{2}
\dot{r}\left(\rho
(\hat{y}(t))\right)\big].
\end{aligned}
\end{equation}
Associated with \eqref{eq2}, we obtain:
\begin{equation}\label{eq3}
u^*(t)=-R^{-1}\big[By(t)+
\frac{1}{2}
\dot{r}\left(\rho
(\hat{y}(t))\right)
+\bar{B}\hat{y}(t)\big].
\end{equation}
Therefore, the optimal state $x^*(\cdot)$ and the corresponding conditional expectation term $\hat{x}(\cdot)$ solve the following equation
\begin{equation}\label{sde11}
\left\{
\begin{aligned}
\mathrm{d} x^*(t)= & \Big\{A x^*(t)+a\left( \hat{x}(t)\right)-R^{-1}B^2
y(t)-\frac{1}{2}R^{-1}(B+\bar{B})
\dot{r}\left(\rho
(\hat{y}(t))\right)\\
&-R^{-1}\bar{B}(2B+\bar{B})\hat{y}(t)
\Big\} \mathrm{d}t
+\sigma \mathrm{d} W(t)+\sigma_0 \mathrm{d} W^0(t), \ t \in(0, T], \\
x^*(0) = &\xi,
\end{aligned}
\right.
\end{equation}
and
\begin{equation}\label{sde13}
\left\{
\begin{aligned}
\mathrm{d}\hat{x}(t)=&\Big\{A \hat{x}(t)+a\left( \hat{x}(t)\right)-R^{-1}(B+\bar{B})^2\hat{y}(t)
-\frac{1}{2}R^{-1}(B+\bar{B})
\dot{r}\left(\rho
(\hat{y}(t))\right)
\Big\} \mathrm{d}t\\
&
+\sigma_0 \mathrm{d} W^0(t), \ t \in(0, T], \\
\hat{x}(0) = &\mathbb{E}^{mf}[\xi].
\end{aligned}
\right.
\end{equation}
Combining \eqref{bsde1}, \eqref{sde11} and \eqref{sde13}, we derive the following stochastic Hamiltonian system where adding $\hat{x}(\cdot)$ as a new state process into the system. This system, which includes $\hat{x}(\cdot)$, is referred to as a mean-field FBSDE (MF-FBSDE):
\begin{equation}\label{fbsde1}
\left\{
\begin{aligned}
\mathrm{d} x^*(t)= & \Big\{A x^*(t)+a\left( \hat{x}(t)\right)-R^{-1}B^2
y(t)-\frac{1}{2}R^{-1}(B+\bar{B})
\dot{r}\left(\rho
(\hat{y}(t))\right)\\
&-R^{-1}\bar{B}(2B+\bar{B})\hat{y}(t)
\Big\} \mathrm{d}t
+\sigma \mathrm{d} W(t)+\sigma_0 \mathrm{d} W^0(t), \ t \in(0, T], \\
\mathrm{d}\hat{x}(t)=&\Big\{A \hat{x}(t)+a\left( \hat{x}(t)\right)-R^{-1}(B+\bar{B})^2\hat{y}(t)
-\frac{1}{2}R^{-1}(B+\bar{B})
\dot{r}\left(\rho
(\hat{y}(t))\right)
\Big\} \mathrm{d}t\\
&
+\sigma_0 \mathrm{d} W^0(t), \ t \in(0, T], \\
\mathrm{d} y(t)= & -\Big\{Ay(t)+\dot{a}\left(\hat{x}(t)\right)
\hat{y}(t) +Q x^*(t)
+ \frac{1}{2}
\dot{q}\left(\hat{x}(t)\right)\Big\}
d t\\
&+ z(t) \mathrm{d} W(t)+ z_0(t) \mathrm{d} W^0(t), \quad t \in[0, T), \\
x^*(0) = &\xi,\
\hat{x}(0) = \mathbb{E}^{mf}[\xi],\
y(T)=  G x^*(T)+\frac{1}{2}
\dot{g}\left( \hat{x}(T)\right).
\end{aligned}
\right.
\end{equation}
We say $(x^*(\cdot),\hat x(\cdot),y(\cdot),z(\cdot),z_0(\cdot))$ is a strong solution to MF-FBSDE \eqref{fbsde1} if 
\begin{equation*}
x^*(\cdot),\hat x(\cdot),y(\cdot)\in L_{\mathbb{F}^{mp}}^2(\Omega^{mf};C(0,T;\mathbb{R}))\quad\text{and}\quad z(\cdot),z_0(\cdot)\in L^2_{\mathbb{F}^{mp}}(0,T
      ;\mathbb{R}),
\end{equation*}
and, $\mathbb P^{mf}$-a.s. and for all $t\in [0,T]$,
\begin{equation*}
\begin{aligned}
x^*(t)=&\xi+\sigma W(t)+\sigma_0 W^0(t)+\int_0^t \Big[Ax^*(s)+a(\hat x(s))-R^{-1}B^2y(s)\\
&-\frac{1}{2}R^{-1}(B+\bar{B})
\dot{r}\left(\rho
(\hat{y}(s))\right)-R^{-1}\bar{B}(2B+\bar{B})\hat{y}(s)
\Big] \mathrm{d}s,\\
\hat x(t)=&\mathbb E^{mf}[\xi]+\sigma_0W^0(t)+\int_0^t\Big[A \hat{x}(s)+a\left( \hat{x}(s)\right)-R^{-1}(B+\bar{B})^2\hat{y}(s)\\
&-\frac{1}{2}R^{-1}B
\dot{r}\left(\rho
(\hat{y}(s))\right)\Big]\mathrm{d} s,\\
y(t)=&G x^*(T)+\frac{1}{2}
\dot{g}\left( \hat{x}(T)\right)-\int_t^Tz(s)\mathrm{d}W(s)-\int_t^Tz_0(s)\mathrm{d}W^0(s)\\
&+\int_t^T\Big[Ay(s)+\dot{a}\left(\hat{x}(s)\right)
\hat{y}(s) +Q x^*(s) +\frac{1}{2}
\dot{q}\left(\hat{x}(s)\right)\Big]
d s.
\end{aligned}
\end{equation*}
Similarly, we can define strong solutions to other MF-FBSDEs.

We conjecture that $$y(t)=P(t)(x^*(t)-\hat{x}(t))+\varphi(t),\ \text{for any $t\in [0,T]$},$$
where $P(\cdot)$ is a deterministic function and $\varphi(\cdot)$ is a BSDE generated by $W^0(\cdot)$ only.
Following the conjecture, we have
$\hat{y}(t)=\varphi(t),\ \text{for any $t\in [0,T]$}.$
Applying It\^o formula, we get:
\begin{equation*}
\begin{aligned}
dy(t)=&\dot{P}(t)(x^*(t)-\hat{x}(t))\mathrm{d}t+P(t)\mathrm{d}(x^*(t)-\hat{x}(t))
+\mathrm{d}\varphi(t)\\
=&\dot{P}(t)(x^*(t)-\hat{x}(t))\mathrm{d}t+P(t)\Big\{A (x^*(t)-\hat{x}(t))-R^{-1}B^2
(y(t)-\hat{y}(t))
\\
&+ P(t)\sigma \mathrm{d}W(t)
+\mathrm{d}\varphi(t)\\
=& -\big[Ay(t)+\dot{a}\left(\hat{x}(t)\right)
\hat{y}(t) +Q x^*(t)
+ \frac{1}{2}
\dot{q}\left(\hat{x}(t)\right)\big]
\mathrm{d}t
 + z(t)\mathrm{d}W(t)+ z_0(t)\mathrm{d}W^0(t).\\
\end{aligned}
\end{equation*}
And then, the above formula equals
\begin{equation*}
\begin{aligned}
&\dot{P}(t)(x^*(t)-\hat{x}(t))\mathrm{d}t+P(t)\Big\{A(x^*(t)-\hat{x}(t))-R^{-1}B^2P(t)(x^*(t)-\hat{x}(t))\Big\} \mathrm{d}t \\
&+ P(t)\sigma \mathrm{d}W(t)
+\mathrm{d}\varphi(t)\\
=& -\Big\{A\big[P(t)(x^*(t)-\hat{x}(t))+\varphi(t)\big]
+\dot{a}\left(\hat{x}(t)\right)
\big[P(t)(x^*(t)-\hat{x}(t))+\varphi(t)\big]\\
&+Q x^*(t)
+ \frac{1}{2}
\dot{q}\left(\hat{x}(t)\right)\Big\}
\mathrm{d} t+ z(t) \mathrm{d} W(t)+ z_0(t) \mathrm{d} W^0(t).
\end{aligned}
\end{equation*}
Therefore, $P(\cdot)$ and $\varphi(\cdot)$ should satisfy the following equations, respectively,
\begin{equation}\label{rica1}
\left\{
\begin{aligned}
&\dot{P}(t)+2AP(t)+Q-R^{-1}B^2P^2(t)
=0,\\
&P(T)=G,
\end{aligned}
\right.
\end{equation}
and
\begin{equation}\label{bsde2}
\left\{
\begin{aligned}
\mathrm{d}\varphi(t)=&-\big[
A\varphi(t)+\dot{a}\left(\hat{x}(t)\right)\varphi(t)+Q\hat{x}(t)+\frac{1}{2}\dot{q}(\hat{x}(t))\big]
\mathrm{d}t+\Lambda_0(t)\mathrm{d}W^0(t),\ t\in[0,T),\\
\varphi(T)=&G\hat{x}(T)+\frac{1}{2}\dot{g}(\hat{x}(T)),
\end{aligned}
\right.
\end{equation}
Substituting  $\hat{y}(\cdot)=\varphi(\cdot)$ into equation \eqref{sde13}, we have $\hat{x}(\cdot)$ satisfies
\begin{equation}\label{sde8}
\left\{
\begin{aligned}
\mathrm{d} \hat{x}(t)= & \big[A\hat{x}(t)
+a\left(\hat{x}(t)\right)
-R^{-1}(B+\bar{B})^2\varphi(t)\\
&
-R^{-1}(B+\bar{B})\tilde{r}\left(t,
\varphi(t)\right)
\big] \mathrm{d}t+\sigma_0 \mathrm{d} W^0(t), \ t \in(0, T], \\
\hat{x}(0)= &\mathbb{E}^{mf}[\xi].
\end{aligned}
\right.
\end{equation}
where $\tilde{r}(x):=\frac{1}{2}\dot{r}\left(\rho\left(x\right)\right)$ and it can be verified that $\tilde{r}$ is a Lipschitz continuous function with respect to $x$.
Combining with \eqref{bsde2}, we derive the following FBSDE:
\begin{equation}\label{fbsde2}
\left\{
\begin{aligned}
\mathrm{d}\hat{x}(t)= & \big[A\hat{x}(t)
+a\left(\hat{x}(t)\right)
-R^{-1}(B+\bar{B})^2\varphi(t)\\
&
-R^{-1}(B+\bar{B})\tilde{r}\left(
\varphi(t)\right)
\big] \mathrm{d}t+\sigma_0\mathrm{d}W^0(t), \ t \in(0, T], \\
\mathrm{d}\varphi(t)=&-\big[
A\varphi(t)+\dot{a}\left(\hat{x}(t)\right)\varphi(t)+Q\hat{x}(t)+\frac{1}{2}\dot{q}(\hat{x}(t))\big]
\mathrm{d}t+\Lambda_0(t)\mathrm{d}W^0(t),\ t\in[0,T),\\
\hat{x}(0)= &\mathbb{E}^{mf}[\xi],\
\varphi(T)=G\hat{x}(T)+
\frac{1}{2}\dot{g}(\hat{x}(T)).
\end{aligned}
\right.
\end{equation}
The system above is a coupled FBSDE system driven only by the common noise $W^0$. It is non-degenerate if $\sigma_0\neq0$. We shall follow the ideas in \cite{delarue2002existence,ma1994four} to show that it is globally well-posed.

\begin{theorem}\label{theo2}
Assume that {\rm (A1)-(A3)} hold. Then FBSDE \eqref{fbsde2} admits a unique strong solution on $[0,T]$. 
\end{theorem}
\noindent{\bf Proof}\quad
Suppose that $(\hat{x}^\eta(\cdot),
\varphi^\eta(\cdot),
\Lambda_0^\eta(\cdot))$ is a strong solution of the following FBSDE:
\begin{equation}\label{fbsde7}
\left\{
\begin{aligned}
\mathrm{d} \hat{x}^\eta(t)= & \Big\{
-R^{-1}(B+\bar{B})^2
\tilde{\eta}_1(t)\varphi^\eta(t)
+
\big[a^\eta\left(t,
\hat{x}^\eta(t)\right)\\
&-R^{-1}(B+\bar{B})\tilde{r}^\eta\left(t,\hat{x}^\eta(t),
\varphi^\eta(t)\right)\big]\tilde{\eta}_1(t)
\Big\} \mathrm{d}t+\tilde{\eta}_1(t)\sigma_0 \mathrm{d} W^0(t), \ t \in(0, T], \\
\mathrm{d}\varphi^\eta(t)=&-\Big\{\big[
A+\tilde{a}^\eta\left(t,
\hat{x}^\eta(t)\right)-R^{-1}(B+\bar{B})^2\tilde{\eta}_2(t)
\big]
\varphi^\eta(t)+\tilde{\eta}_2(t)\big[a^\eta(t,\hat{x}^{\eta}(t))\\
&-R^{-1}(B+\bar{B})\tilde{r}^\eta\left(t,\hat{x}^\eta(t),
\varphi^\eta(t)\right)\big]
+\tilde{q}^\eta\left(t,
\hat{x}^\eta(t)\right)
\Big\}\mathrm{d}t\\
&+\Lambda_0^\eta(t)\mathrm{d}W^0(t),\ t\in[0,T),\\
\hat{x}^\eta(0)= &\mathbb{E}^{mf}[\xi],\
\varphi^\eta(T)=\tilde{g}^\eta
\left(T,\hat{x}^\eta(T)\right).
\end{aligned}
\right.
\end{equation}
Here,
\begin{equation}\label{eq11}
\begin{aligned}
a^\eta\left(t,x\right):=&a\left(
\tilde{\eta}_1^{-1}(t)x\right),\
\tilde{a}^\eta\left(t,x\right):=\dot{a}
\left(
\tilde{\eta}^{-1}_1(t)x\right),\\
\tilde{r}^\eta(t,x,y):=&
\tilde{r}\left(y+\eta_2(t)\eta_1^{-1}(t)x\right),\
\tilde{q}^\eta\left(t,x\right)
:=\frac{1}{2}
\dot{q}\left(
\tilde{\eta}_1^{-1}(t)x\right),\
\tilde{g}^\eta\left(t,x\right)
:=\frac{1}{2}
\dot{g}\left(
\tilde{\eta}_1^{-1}(t)x\right),
\end{aligned}
\end{equation}
where $\tilde{\eta}_1(\cdot)$ and $\tilde{\eta}_2(\cdot)$ are given by the following equations:
\begin{equation}\label{ode1}
\left\{
\begin{aligned}
\mathrm{d} \tilde{\eta}_1(t)= & -\big[
A-R^{-1}(B+\bar{B})^2\tilde{\eta}_2(t)
\big]
\tilde{\eta}_1(t)\mathrm{d}t,\\
\tilde{\eta}_1(T)= &1,
\end{aligned}
\right.
\end{equation}
and
\begin{equation}\label{ode2}
\left\{
\begin{aligned}
\mathrm{d} \tilde{\eta}_2(t)= & -\Big\{
\big[2A
+\tilde{a}^{\eta}\left(t,\hat{x}^{\eta}(t)\right)\big]
\tilde{\eta}_2(t)-R^{-1}(B+\bar{B})^2
\tilde{\eta}^2_2(t)+Q
\Big\}\mathrm{d}t,\\
\tilde{\eta}_2(T)= &G.
\end{aligned}
\right.
\end{equation}
Note that \eqref{fbsde7}-\eqref{ode2} do not involve $(\hat x(\cdot),\varphi(\cdot),\Lambda_0(\cdot))$. Under (A1) and (A3), it can be verified that $a^\eta$, $\tilde{a}^\eta$, $\tilde{q}^\eta$, $\tilde{g}^\eta$ are bounded uniformly Lipschitz continuous with respect to $x$, and $\tilde{r}^\eta$ are bounded uniformly Lipschitz continuous with respect to $x$ and $y$. 
Let 
\begin{equation}\label{eq14}
\tilde{\hat{x}}(\cdot):=
\tilde{\eta}_1^{-1}(\cdot)
\hat{x}^\eta(\cdot),\
\tilde{\varphi}(\cdot):=\varphi^{\eta}
(\cdot)+\tilde{\eta}_2(\cdot)
\tilde{\eta}_1^{-1}(\cdot)
\hat{x}^\eta(\cdot),\
\tilde{\Lambda}_0(\cdot):=
\Lambda_0^\eta(\cdot)
+\tilde{\eta}_2(\cdot)\sigma_0.
\end{equation}
Then, we can verify that $(\tilde{\hat{x}}(\cdot),\tilde{\varphi}(\cdot),
\tilde{\Lambda}_0(\cdot))$ solves FBSDE \eqref{fbsde2}.

Suppose that FBSDE system \eqref{fbsde2} admits a strong solution $(\hat{x}(\cdot),
\varphi(\cdot),
\Lambda_0(\cdot))$. Introduce the following auxiliary equations:
\begin{equation}\label{sde15}
\left\{
\begin{aligned}
\mathrm{d} \eta_1(t)= & -\big[
A-R^{-1}(B+\bar{B})^2\eta_2(t)
\big]
\eta_1(t)\mathrm{d}t,\\
\eta_1(T)= &1,
\end{aligned}
\right.
\end{equation}
and
\begin{equation}\label{sde16}
\left\{
\begin{aligned}
\mathrm{d}\eta_2(t)= & -\Big\{
\big[2A
+\dot{a}\left(\hat{x}(t)\right)\big]
\eta_2(t)-R^{-1}(B+\bar{B})^2
\eta^2_2(t)+Q
\Big\}\mathrm{d}t,\\
\eta_2(T)= &G.
\end{aligned}
\right.
\end{equation}
Under (A1) and (A3), fix $\omega^0\in\Omega^0$, using \cite[Chapter 6, Theorem 7.2]{yong1999stochastic}, we have equation \eqref{sde16} admits a unique solution $\eta_2(\cdot)(\omega^0)\in C([0,T];\mathbb{R}_{+})$. And then, by solving the linear backward ODE \eqref{sde15}, there exists a constant $M>0$ such that $M\leq\eta_1(\cdot)\leq\frac{1}{M}$.
Consider the transformation: 
\begin{equation*}
\tilde{\hat{x}}^\eta(\cdot)
:=\eta_1(\cdot)\hat{x}(\cdot),\
\tilde{\varphi}^\eta(\cdot):=\varphi(\cdot)
-\eta_2(\cdot)\hat{x}(\cdot),\
\tilde{\Lambda}_0^\eta(\cdot):=\Lambda_0(\cdot)
-\eta_2(\cdot)\sigma_0.
\end{equation*}
By straightforward calculation, we can check that $(\tilde{\hat{x}}^\eta(\cdot),\tilde{\varphi}^\eta(\cdot),
\tilde{\Lambda}_0^\eta(\cdot))$ indeed satisfies FBSDE \eqref{fbsde7} with $(\tilde \eta_1,\tilde\eta_2)$ replaced by $(\eta_1,\eta_2)$. Finally, we can verify that $\eta_2$ satisfies
\begin{equation*}
\left\{
\begin{aligned}
\mathrm{d} \eta_2(t)= & -\Big\{
\big[2A
+\dot{a}\left(\eta_1^{-1}(t)\tilde{\hat{x}}^{\eta}(t)\right)\big]
\tilde{\eta}_2(t)-R^{-1}(B+\bar{B})^2
\tilde{\eta}^2_2(t)+Q
\Big\}\mathrm{d}t,\\
\tilde{\eta}_2(T)= &G.
\end{aligned}
\right.
\end{equation*}
which yields that the well-posedness between FBSDEs \eqref{fbsde2} and \eqref{fbsde7} are equivalent. 

Noting that under (A1)-(A3), the drift term of the SDE and the generator of the BSDE in \eqref{fbsde7} are bounded with respect to $\hat{x}^{\eta}$ and exhibit linear growth with respect to $\varphi^{\eta}$, while the terminal condition of the BSDE is bounded, and the diffusion term of the SDE is non-degenerate. Therefore, according to \cite[Theorem 2.6]{delarue2002existence}, it follows that FBSDE \eqref{fbsde7} admits a unique strong solution $(\hat{x}^\eta(\cdot),
\varphi^\eta(\cdot),
\Lambda_0^\eta(\cdot))$. And thus, by the equivalency between FBSDEs \eqref{fbsde2} and \eqref{fbsde7}, we have the existence and uniqueness of strong solutions to FBSDE \eqref{fbsde2}.

\hfill$\Box$
\begin{remark}

Consider a special case where $\dot{a}=\bar{B}=\dot{q}=\dot{r}
=\dot{g} =0$. In this situation, it can be checked that $\varphi(\cdot)$ becomes a linear backward ODE, and Problem (MF) degenerates into a classic LQ control problem in \cite[Chapter 6]{yong1999stochastic}.
\end{remark}
\begin{theorem}\label{theo3}
Assume that {\rm (A1)-(A3)} hold. Then MF-FBSDE \eqref{fbsde1} admits a unique strong solution.
\end{theorem}
\noindent{\bf Proof}\quad
\emph{Existence.} Let  $(\hat{x}(\cdot),\varphi(\cdot))$ be the unique solution to FBSDE \eqref{fbsde2}. Consider the following equation
\begin{equation}\label{sde12}
\left\{
\begin{aligned}
\mathrm{d}x^{\hat{x},\varphi}(t)= &\Big\{A x^{\hat{x},\varphi}(t)+a\left( \hat{x}(t)\right)-R^{-1}B^2
P(t)(x^{\hat{x},\varphi}(t)-\hat{x}(t))\\
&-R^{-1}(B+\bar{B})
\tilde{r}\left(\varphi(t)\right)-R^{-1}(B+\bar{B})^2\varphi(t)
\Big\} \mathrm{d}t\\
&+\sigma \mathrm{d} W(t)+\sigma_0 \mathrm{d} W^0(t), \ t \in(0, T], \\
x^{\hat{x},\varphi}(0) = &\xi,
\end{aligned}
\right.
\end{equation}
which is well-posed on $[0,T]$. Taking the conditional expectation of \eqref{sde12} on $\mathcal{F}^{W^0}_t$, for $t\in[0, T]$, we have:
\begin{equation*}
\left\{
\begin{aligned}
\mathrm{d}\mathbb{E}^{mf}[
x^{\hat{x},\varphi}(t)\vert
\mathcal{F}^{W^0}_t]= & \Big\{A \mathbb{E}^{mf}[
x^{\hat{x},\varphi}(t)\vert
\mathcal{F}^{W^0}_t]+a\left( \hat{x}(t)\right)-R^{-1}(B+\bar{B})^2
\varphi(t)\\
&
-R^{-1}(B+\bar{B})\tilde{r}\left(\varphi(t)
\right)
\Big\} \mathrm{d}t +\sigma_0 \mathrm{d}W^0(t), \ t \in(0, T], \\
\mathbb{E}^{mf}[
x^{\hat{x},\varphi}
(0)] = &\mathbb{E}^{mf}[\xi].
\end{aligned}
\right.
\end{equation*}
Comparing with \eqref{sde8}, it yields that
\begin{equation*}
\left\{
\begin{aligned}
\mathrm{d}\left(\mathbb{E}^{mf}[
x^{\hat{x},\varphi}(t)\vert
\mathcal{F}^{W^0}_t]-\hat{x}(t)\right)= & A
 \big(\mathbb{E}^{mf}[
x^{\hat{x},\varphi}(t)\vert
\mathcal{F}^{W^0}_t]-\hat{x}(t)\big)
\mathrm{d}t, \ t \in(0, T], \\
\mathbb{E}^{mf}[
x^{\hat{x},\varphi}
(0)]-\hat{x}(0) = &0.
\end{aligned}
\right.
\end{equation*}
Hence, we have $\hat{x}(t)=
\mathbb{E}^{mf}[x^{\hat{x},\varphi}(t)
\vert\mathcal{F}^{W^0}_t]$ for $ t\in[0,T]$.
Define
$$y^{\hat{x},\varphi}(\cdot):=
P(\cdot)(x^{\hat{x},\varphi}(\cdot)-\hat{x}(\cdot))
+\varphi(\cdot),\
z^{\hat{x},\varphi}(\cdot):=
P(\cdot)\sigma,\
z^{\hat{x},\varphi}_0(\cdot):=
P(\cdot)\sigma_0+\Lambda_0(\cdot).
$$
Applying It\^o's formula to $y^{\hat{x},\varphi}(\cdot)$, it follows that $\left(x^{\hat{x},\varphi}(\cdot),
\hat{x}(\cdot),
y^{\hat{x},\varphi}(\cdot),
z^{\hat{x},\varphi}(\cdot),
z^{\hat{x},\varphi}_0(\cdot)\right)$ satisfies FBSDE \eqref{fbsde1}.

\emph{Uniqueness.} Assume that $(\bar{x}^*(\cdot),
\bar{\hat{x}}(\cdot),
\bar{y}(\cdot),
\bar{z}(\cdot),
\bar{z}_0(\cdot))$ is another solution to FBSDE \eqref{fbsde1}.
Similarly, we get:
\begin{equation*}
\left\{
\begin{aligned}
\mathrm{d}\left(\mathbb{E}^{mf}[
\bar{x}^*(t)\vert
\mathcal{F}^{W^0}_t]-\bar{\hat{x}}(t)
\right)= & A
\big(\mathbb{E}^{mf}[
\bar{x}^*(t)\vert
\mathcal{F}^{W^0}_t]-\bar{\hat{x}}(t)
\big)
\mathrm{d}t, \ t \in(0, T], \\
\mathbb{E}^{mf}[
\bar{x}^*(0)]-\bar{\hat{x}}(0)
 = &0,
\end{aligned}
\right.
\end{equation*}
which means that $\bar{\hat{x}}(t)=\mathbb{E}^{mf}[
\bar{x}^*(t)\vert
\mathcal{F}^{W^0}_t]$ for $t\in[0,T]$.
Let $\bar{y}(\cdot):=\bar{P}(\cdot)(
\bar{x}^*(\cdot)-\bar{\hat{x}}(\cdot))
+\bar{\varphi}(\cdot)$.
And then, we have $\mathbb{E}^{mf}[
\bar{y}(t)\vert
\mathcal{F}^{W^0}_t]=\bar{\varphi}(t)$, for any $t\in[0,T]$. Due to the previous derivation, we have $\bar{P}(\cdot)$ and $(\bar{\hat{x}}(\cdot),
\bar{\varphi}(\cdot))$ satisfy \eqref{rica1} and \eqref{fbsde2}, respectively. Therefore, the desired
uniqueness follows from that of the solutions to \eqref{rica1} and \eqref{fbsde2}.
\hfill$\Box$

To summarize the above discussion, we have the following theorem.
\begin{theorem}\label{thm:mfu}
Assume that {\rm (A1)-(A3)} hold. Let $P(\cdot)$ and $(\hat{x}(\cdot),\varphi(\cdot))$ be the unique solution to \eqref{rica1} and \eqref{fbsde2}, respectively. Then the optimal control for Problem (MF) has the following form for any $t\in[0,T]$,
\begin{equation}\label{u1}
\begin{aligned}
u^*(t)=&-R^{-1}BP(t)(x^*(t)
-\hat{x}(t))+\rho(\varphi(t))\\
=&-R^{-1}\big[BP(t)(x^*(t)-\hat{x}(t))
+(B+\bar{B})\varphi(t)+\tilde{r}\left(
\varphi(t)\right)\big],
\end{aligned}
\end{equation}
where the optimal state $x^*(\cdot)$ solves
\begin{equation}\label{sde2}
\left\{
\begin{aligned}
\mathrm{d}x^*(t)= & \big\{A x^*(t)+a\left( \hat{x}(t)\right)-R^{-1}B^2P(t)
(x^*(t)-\hat{x}(t))+(B+\bar{B})\rho\left(\varphi(t)\right)
\big\} \mathrm{d}t\\
&+\sigma \mathrm{d} W(t) +\sigma_0 \mathrm{d} W^0(t), \ t \in(0, T], \\
x^*(0) = &\xi.
\end{aligned}
\right.
\end{equation}
\end{theorem}
\begin{example}\label{example1}
Consider the following system
\begin{equation*}
\left\{
\begin{aligned}
   \mathrm{d}x(t)=&u(t)dt+\sigma \mathrm{d}W(t)+\sigma_0 \mathrm{d}W^0(t),\ t\in(0,T],\\
    x(0)=&\xi,
\end{aligned}
\right.
\end{equation*}
and the cost functional being
\begin{equation*}
    J(\mu;u(\cdot))=\mathbb{E}^{mf}\left[\int_0^T\frac{1}{2}u^2(t)\mathrm{d}t-\big\vert\mathbb{E}^{mf}[x(T)\vert\mathcal{F}_T^0]\big\vert^2\right],
\end{equation*}
where $\sigma_0>0$, $\xi\sim \mu$ and $\mathbb{E}^{mf}[\xi]=x_0$.
For any $l>0$, we take $\mathbb{E}^{mf}[u_l(t)\vert\mathcal{F}_t^0]=l$, $t\in[0,T]$.
Let $x_l(\cdot)$ be the corresponding trajectory. Then
\begin{equation*}
  \mathbb{E}^{mf}\big\vert\mathbb{E}^{mf}[x_l(T)\vert\mathcal{F}_T^0]\big\vert^2=\sigma_0^2T+(lT+x_0)^2.
\end{equation*}
Thus,
\begin{equation*}
\begin{aligned}
    J(\mu;u_l(\cdot))=&\mathbb{E}^{mf}\left\{\int_0^T\left[\frac{1}{2}(u_l(t)-\mathbb{E}^{mf}[u_l(t)\vert\mathcal{F}_t^0])^2+\frac{1}{2}\big\vert\mathbb{E}^{mf}[u_l(t)\vert\mathcal{F}_t^0]\big\vert^2\right]\mathrm{d}t-\big\vert\mathbb{E}^{mf}[x_l(T)\vert\mathcal{F}_T^0]\big\vert^2\right\}\\
    =&\mathbb{E}^{mf}\left[\int_0^T\frac{1}{2}(u_l(t)-l)^2\mathrm{d}t\right]+\frac{1}{2}l^2T-\sigma_0^2T-(lT+x_0)^2,
\end{aligned}    
\end{equation*}
which yields that $V(0,x_0)\rightarrow -\infty$, as $l\rightarrow \infty$.
\end{example}
The above example demonstrates that if the cost functional is strictly concave with respect to the conditional expectation term (i.e., its second derivative is less than or equal to a negative constant), an optimal control may not exist. However, the following example shows that under our results, even if the cost functional is concave with respect to the conditional expectation term, an optimal control still exists and is unique.
\begin{example}\label{example2}
Consider control problem \eqref{sde1}-\eqref{eq:valueint} with $Q=G=0,\ R=2,\ a(\cdot),q(\cdot),r(\cdot),g(\cdot)$ satisfy $f(x)=-\frac{1}{2}x^2$ for $\vert x\vert\leq 1$, $\vert f'(x)\vert=0$ for $\vert x\vert\geq 2$ and $\vert f'(x)\vert\leq 1$, for $x\in \mathbb{R}$, where $f(\cdot)=a(\cdot),q(\cdot),r(\cdot),g(\cdot)$. It can be verified that this situation satisfies {\rm (A1)-(A3)}. In this case, functions $a(\cdot)$, $q(\cdot)$, $r(\cdot)$ and $g(\cdot)$ are non-convex. Based on the previous analysis, the above control problem has a unique optimal control $u^*(\cdot)$.
\end{example}
\subsection{LQ case}
Now, let us consider the LQ case where $a(x)=\bar{A}x$, $q(x)=\bar{Q}x^2$, $r(x)=\bar{R}x^2$ and $g(x)=\bar{G}x^2$. In this case, the state equation and the cost functional take the following forms, respectively.
\begin{equation*}
\left\{\begin{aligned}
\mathrm{d}x(t)= & \left\{A x(t)+\bar{A} \mathbb{E}[x(t)|\mathcal{F}^{W^0}_t]+B u(t)+ \bar{B} \mathbb{E}[u(t)|\mathcal{F}^{W^0}_t]
\right\} \mathrm{d}t +\sigma d W(t) +\sigma_0 d W^0(t), \ t \in(0, T], \\
x(0)= & \xi,
\end{aligned}\right.
\end{equation*}
and
\begin{equation*}
\begin{aligned}
J(\xi ; u(\cdot))= & \mathbb{E}^{mf}\Big\{\int_0^T\left[Q x^2(t)
+\bar{Q}\left(\mathbb{E}^{mf}[x(t)\vert
\mathcal{F}^{W^0}_t]\right)^2
+ R u^2(t)
+\bar{R}\left( \mathbb{E}^{mf}[u(t)\vert
\mathcal{F}^{W^0}_t] \right)^2\right] \mathrm{d} t  \\
& + G x^2(T)+ \bar{G}\left( \mathbb{E}^{mf}[x(T)\vert\mathcal{F}^{W^0}_t] \right)^2\Big\}.
\end{aligned}
\end{equation*}
Assume that $Q+\bar{Q}\geq0$, $R+\bar{R}>0$, $R>0$ and $G+\bar{G}\geq0$. 
And then, we have the optimal control satisfies the following form:
$$u^*(t)=-R^{-1}B(y(t)-\hat{y}(t))
-(R+\bar{R})^{-1}(B+\bar{B})\hat{y}
(t),$$
where $y(\cdot)$ reduces to
\begin{equation*}
\left\{
\begin{aligned}
\mathrm{d} y(t)= & -[Ay(t)+\bar{A}
\hat{y}(t) +Q x^*(t)+ \bar{Q}\hat{x}(t)
]\mathrm{d}t
+ z(t)\mathrm{d}W(t)+ z_0(t)\mathrm{d}W^0(t), \ t \in[0, T), \\
y(T)= & G x^*(T)+\bar{G} \hat{x}(T).
\end{aligned}
\right.
\end{equation*}
Furthermore, FBSDE \eqref{fbsde2}
becomes
\begin{equation*}
\left\{
\begin{aligned}
\mathrm{d}\hat{x}(t)= & \big\{(A+\bar{A})\hat{x}(t)
-(R+\bar{R})^{-1}(B+\bar{B})^2
\varphi(t)
\big\} \mathrm{d}t
+\sigma_0 \mathrm{d} W^0(t), \ t \in(0, T], \\
\mathrm{d}\varphi(t)=&-\big\{
(A+\bar{A})\varphi(t)+(Q+\bar{Q})\hat{x}(t)
\big\}\mathrm{d}t+\Lambda_0(t)\mathrm{d}W^0(t),\ t\in[0,T),\\
\hat{x}(0)= &\mathbb{E}[\xi],\
\varphi(T)=(G+\bar{G})\hat{x}(T).
\end{aligned}
\right.
\end{equation*}
It can be verified that $$\varphi(\cdot)=\eta(\cdot)\hat{x}(\cdot)+\phi(\cdot),\
\Lambda_0(\cdot)=\eta(\cdot)\sigma_0,$$
where $\eta(\cdot)$ and $\phi(\cdot)$ solve the following equations, respectively.
\begin{equation*}
\left\{
\begin{aligned}
&\mathrm{d}\eta(t)=-\Big\{2(A+\bar{A})
\eta(t)+(Q+\bar{Q})-(R+\bar{R})^{-1}
(B+\bar{B})^2\eta^2(t)\Big\}\mathrm{d}t,\\
&\eta(T)=G+\bar{G},
\end{aligned}
\right.
\end{equation*}
and
\begin{equation*}
\left\{
\begin{aligned}
\mathrm{d}\phi(t)=&-\Big\{
A+\bar{A}
-(R+\bar{R})^{-1}(B+\bar{B})^2
\eta(t)\Big\}\phi(t)\mathrm{d}t,\ t\in[0,T),\\
\phi(T)=&0.
\end{aligned}
\right.
\end{equation*}
From \cite{yong1999stochastic}, we have $\eta(\cdot)$ admits a unique solution. Therefore,
$$y(\cdot)=P(\cdot)(x(\cdot)-
\hat{x}(\cdot))+\eta(\cdot)\hat{x}(\cdot)
+\phi(\cdot).$$
Noting that $\phi(\cdot)$ is no longer dependent on $\hat{x}(\cdot)$, which means that we have completely decoupled $y(\cdot)$ used by $x^*(\cdot)$ and $\hat x(\cdot)$.

Specifically, if there is no common noise in the system, the corresponding conclusion becomes the result in \cite{yong2013linear}.
\subsection{Hamilton-Jacobi equation}\label{sec1}
In this section, we will deal with the well-posedness of the Hamilton-Jacobi equation \eqref{pde3} corresponding to Problem (MF).

Note that the Hamilton-Jacobi equation \eqref{pde3} is an infinite-dimensional PDE and its solution $V$ is a map: $[0,T]\times\mathcal{P}_2(\mathbb R)\to \mathbb R$. To study the well-posedness of the Hamilton-Jacobi equation \eqref{pde3}, we need the well-posedness of the following two finite dimensional PDEs which can be regarded as the decoupling fields of FBSDEs
\eqref{fbsde1} and \eqref{fbsde2}, respectively. For the concept of the decoupling field, please refer to \cite[Chapter 8, Definition 8.3.2]{zhang2017backward}.
\begin{equation}\label{pde2}
\left\{
\begin{aligned}
&\partial_tU(t,x,\hat{x})+\partial_x
U(t,x,\hat{x})\Big[Ax+a\left(\hat{x}
\right)
-\frac{1}{2}R^{-1}B^2(U(t,x,\hat x)-U(t,\hat x,\hat x))\\
&+(B+\bar{B})\rho\left(\frac{1}{2}U(t,\hat{x},\hat{x})
\right)\Big]
+\partial_{\hat{x}}U(t,x,\hat{x})
\Big[A\hat{x}
+a(\hat{x})+(B+\bar{B})\rho\left(\frac{1}{2}U(t,\hat{x},\hat{x})
\right)\Big]\\
&+\frac{\sigma^2_0}{2}\partial_{\hat{x}\hat{x}}
U(t,x,\hat{x})+\sigma^2_0\partial_{x\hat{x}}
U(t,x,\hat{x})+\frac{\sigma^2+\sigma_0^2}{2}\partial_{xx}U(t,x,\hat x)
\\
&+2Qx+\dot{q}(\hat{x})+AU(t,x,\hat{x})
+\dot{a}\left(\hat{x}\right)
U(t,\hat{x},\hat{x})=0, \ t\in [0,T),\\
&U(T,x,\hat{x})=2Gx+\dot{g}
(\hat{x}),
\end{aligned}
\right.
\end{equation}
and
\begin{equation}\label{pde1}
\left\{
\begin{aligned}
&\partial_t \Phi(t, \hat{x})+\partial_{\hat{x}} \Phi(t, \hat{x})\big[A \hat{x}+a\left(\hat{x}
\right)+(B+\bar{B}) \rho\left(\Phi(t,\hat{x})
\right)\big] \\
&+\frac{\sigma_0^2}{2} \partial_{\hat{x}\hat{x}}
\Phi(t, \hat{x})+\frac{1}{2}\dot{q}\left(
\hat{x}\right) +(A+\dot{a}\left(\hat{x}\right)) \Phi(t, \hat{x})+Q\hat{x}=0, \ t\in [0,T),\\
&\Phi(T,\hat{x})=G\hat{x}+\frac{1}{2}\dot{g}(\hat{x}).
\end{aligned}
\right.
\end{equation}
\begin{remark}
We emphasize that PDE \eqref{pde2} is not a classical PDE since it involves the term $U(t,\hat x,\hat x)$. Moreover, it is a degenerate parabolic PDE if $\sigma=0$ (which is a case we consider in the current paper.)
\end{remark}
\begin{remark}\label{remarkinfity}
It should be noted that the Hamilton-Jacobi equation considered in this paper is essentially infinite-dimensional, even in the following simplified case: the state equation satisfies  
$$
\mathrm{d} x(s)=u(s) \mathrm{d} s+\sigma \mathrm{d} W_s+\sigma_0 \mathrm{d} W_s^0, \ x(t)\vert\mathcal{F}_t^0 \sim \mu,\ 0\leq t\leq s\leq T,
$$
and the cost functional $J$ is defined by: 
\begin{equation*}
J(\mu;u(\cdot)):= \mathbb{E}^{mf}\left[\int_t^T\frac{1}{2}\vert u(s)\vert^2\mathrm{d}s+\vert x(T)\vert^2+g(\mathbb{E}^{mf}\left[x(T)\vert\mathcal{F}_{T}^0\right])\Big\vert\mathcal{F}_t^0\right].
\end{equation*} 
We introduce two PDEs:
\begin{equation*}
\left\{
\begin{aligned}
&\dot{p}(t)-\frac{1}{2}p^2(t)=0,\ t\in(0,T],\\
&p(T)=1,
\end{aligned}
\right.
\end{equation*}
and
\begin{equation*}
\left\{
\begin{aligned}
&-\partial_tv(t,\hat{x})-\frac{\sigma_0^2}{2}\partial_{\hat{x}\hat{x}}v(t,\hat{x})+\frac{1}{2}\vert\partial_{\hat{x}}v(t,\hat{x})\vert^2=0,\ t\in(0,T],\\
&v(T,\hat{x})=g(\hat{x}).
\end{aligned}
\right.
\end{equation*}
Therefore, 
\begin{align}
J=&\mathbb{E}^{mf}\left[\int_t^T\frac{1}{2}\vert u(s)\vert^2\mathrm{d}s+p(T)\vert x(T)\vert^2+v(T,\mathbb{E}^{mf}\left[x(T)\vert\mathcal{F}_{T}^0\right])\Big\vert \mathcal{F}_t^0\right]\notag\\
=&\mathbb{E}^{mf}\Big[\int_t^T\big\{\frac{1}{2}\left(u(s)+2p(s)x(s)+\partial_{\hat{x}}v(s,\mathbb{E}^{mf}\left[x(s)\vert\mathcal{F}_{s}^0\right])\right)^2\notag\\
&-2p(s)x(s)\partial_{\hat{x}}v(s,\mathbb{E}^{mf}\left[x(s)\vert\mathcal{F}_{s}^0\right])\big\}\mathrm{d}s\Big\vert\mathcal{F}_t^0\Big]
+p(t)x^2+v(t,\mathbb{E}^{mf}[\xi])+(\sigma^2+\sigma^2_0)\int_t^Tp(s)\mathrm{d}s\notag\\
\geq &\mathbb{E}^{mf}\left[\int_t^T-2p(s)x(s)\partial_{\hat{x}}v(s,\mathbb{E}^{mf}\left[x(s)\vert\mathcal{F}_{s}^0\right])\mathrm{d}s\Big\vert\mathcal{F}_t^0\right]+p(t)x^2+v(t,\mathbb{E}^{mf}[\xi])\notag\\
&+(\sigma^2+\sigma^2_0)\int_t^Tp(s)\mathrm{d}s,\notag
\end{align}
By selecting the optimal control $u^*(t)=-2p(t)x^{*}(t)-\partial_{\hat{x}}v(t,\hat{x}(t))$,  the value function \begin{equation*}
\begin{aligned} 
V:=&\inf_{\alpha}J\\
=&\int_t^T-2p(s)\hat{x}(s)\partial_{\hat{x}}v(s,\hat{x}(s)) \mathrm{d}s+p(t)x^2+v(t,\mathbb{E}^{mf}[\xi])+(\sigma^2+\sigma^2_0)\int_t^Tp(s)ds,
\end{aligned}
\end{equation*}
Due to the presence of the term $\int_t^T-2p(s)\hat{x}(s)\partial_{\hat{x}}v(s,\hat{x}(s))\mathrm{d}s$, the value function $V$ cannot be expressed as a function of the initial state $x$ and the expectation $\mathbb{E}^{mf}[\xi]$. This indicates that the value function $V$ is essentially infinite-dimensional.
\end{remark}
In order to derive the well-posedness of PDE \eqref{pde1}, recalling the proof of Theorem \ref{theo2}, we introduce the following FBSDE for any $(t_0,x_0)\in[0,T]\times\mathbb{R}$:
\begin{equation}\label{fbsde8}
\left\{
\begin{aligned}
\mathrm{d} \hat{x}^{\eta,t_0,x_0}(t)= & \Big\{
-R^{-1}(B+\bar{B})^2
\tilde{\eta}^{t_0,x_0}_1(t)\varphi^{\eta,t_0,x_0}(t)
+
\big[a^\eta\left(t,
\hat{x}^{\eta,t_0,x_0}\right)\\
&-R^{-1}(B+\bar{B})\tilde{r}^\eta\left(t,\hat{x}^{\eta,t_0,x_0},
\varphi^{\eta,t_0,x_0}(t)\right)\big]\tilde{\eta}^{t_0,x_0}_1(t)
\Big\} \mathrm{d}t\\
&+\tilde{\eta}^{t_0,x_0}_1(t)\sigma_0 \mathrm{d} W^0(t), \ t \in(0, T], \\
\mathrm{d}\varphi^{\eta,t_0,x_0}(t)=&-\Big\{\big[
A+\tilde{a}^\eta\left(t,
\hat{x}^{\eta,t_0,x_0}\right)-R^{-1}(B+\bar{B})^2\tilde{\eta}^{t_0,x_0}_2(t)
\big]
\varphi^{\eta,t_0,x_0}(t)\\
&+\tilde{\eta}^{t_0,x_0}_2(t)\big[a^\eta(t,\hat{x}^{\eta,t_0,x_0})
-R^{-1}(B+\bar{B})\tilde{r}^\eta\left(t,\hat{x}^{\eta,t_0,x_0},
\varphi^{\eta,t_0,x_0}(t)\right)\big]\\
&
+\tilde{q}^\eta\left(t,
\hat{x}^{\eta,t_0,x_0}\right)
\Big\}\mathrm{d}t
+\Lambda^{\eta,t_0,x_0}(t)\mathrm{d}W^0(t),\ t\in[0,T),\\
\hat{x}^{\eta,t_0,x_0}(t_0)= &x_0,\
\varphi^{\eta,t_0,x_0}(T)
=\tilde{g}^\eta
\left(T,
\hat{x}^{\eta,t_0,x_0}(T)\right).
\end{aligned}
\right.
\end{equation}
Here, $\tilde{\eta}_1^{t_0,x_0}(\cdot)$ and $\tilde{\eta}_2^{t_0,x_0}(\cdot)$ are the solutions to equations \eqref{sde15} and \eqref{sde16} respectively with $(\hat{x}(\cdot),\varphi(\cdot))$ replaced by
$\left((\tilde{\eta}_1^{t_0,x_0}(\cdot))^{-1}
\hat{x}^{\eta,t_0,x_0}(\cdot)\right.$,
$\left.\varphi^{\eta,t_0,x_0}(\cdot)
+\tilde{\eta}_2^{t_0,x_0}(\cdot)
(\tilde{\eta}_1^{t_0,x_0}(\cdot))^{-1}
\hat{x}^{\eta,t_0,x_0}(\cdot)\right)$, and $a^\eta$, $\tilde{a}^\eta$, $\tilde{q}^\eta$, $\tilde{g}^\eta$, $\tilde{r}^\eta$ are defined by \eqref{eq11} where $\tilde{\eta}_1(\cdot)$ and
$\tilde{\eta}_2(\cdot)$ are replaced by $\tilde{\eta}_1^{t_0,x_0}(\cdot)$ and
$\tilde{\eta}_2^{t_0,x_0}(\cdot)$, respectively. In addition, FBSDE \eqref{fbsde8} is associated to the following PDE:
\begin{equation}\label{pde7}
\left\{
\begin{aligned}
&\partial_t \tilde{\Phi}(t, \hat{x})+\partial_{\hat{x}} \tilde{\Phi}(t, \hat{x})\Big\{
-R^{-1}(B+\bar{B})^2\tilde{\eta}_1(t)\tilde{\Phi}(t, \hat{x})
+\big[a^\eta\left(t,
\hat{x}\right)\\
&-R^{-1}(B+\bar{B})\tilde{r}^\eta\left(t,
\hat{x},
\tilde{\Phi}(t, \hat{x})\right)
\big]\tilde{\eta}_1(t)\Big\} +\frac{(\tilde{\eta}_1(t))^2 \sigma_0^2}{2} \partial_{\hat{x}\hat{x}}
\tilde{\Phi}(t, \hat{x})\\
&+\big[
A+\tilde{a}^\eta\left(t,\hat{x}
\right)-R^{-1}(B+\bar{B})^2\tilde{\eta}_2^{t_0,x_0}(t)\big]
\tilde{\Phi}(t, \hat{x})\\
&+\tilde{\eta}_2(t)
\big[a^\eta\left(t,
\hat{x}\right)
-R^{-1}(B+\bar{B})\tilde{r}^\eta
\left(t,\hat{x},
\tilde{\Phi}(t, \hat{x})\right)
\big]+\tilde{q}^\eta\left(t,
\hat{x}\right)
=0, \\
&\tilde\Phi(T,\hat{x})=\tilde{g}^\eta
(T,\hat{x}).
\end{aligned}
\right.
\end{equation}
Here, $\tilde{\eta}_1(\cdot)$ and $\tilde{\eta}_2(\cdot)$ are the solutions to equations \eqref{sde15} and \eqref{sde16} respectively with $(\hat{x}(\cdot),\varphi(\cdot))$ replaced by
$\left((\tilde{\eta}_1^{-1}(\cdot))\hat x,
\tilde \Phi(\cdot,\hat x)+\tilde{\eta}_2(\cdot)\tilde{\eta}_1^{-1}(\cdot)\hat x\right)$.

Theorem \ref{thm:tildephi} gives the well-posedness of FBSDE \eqref{fbsde8} and the corresponding PDE \eqref{pde7}. The proof follows the similar arguments as in \cite[Theorem 3]{li2023linear}. For completeness, we provide a proof in Appendix.
\begin{theorem}\label{thm:tildephi}
Assume that {\rm (A1)-(A3)} hold. Then FBSDE \eqref{fbsde8} admits a unique strong solution. Define $\tilde{\Phi}(t_0,x_0)
:=\varphi^{\eta,t_0,x_0}(t_0)
$ for any $(t_0,x_0)\in [0,T]\times\mathbb R$. Then we further have the PDE \eqref{pde7} admits a unique classical solution $\tilde{\Phi}\in C^{1,2}([0,T]\times \mathbb{R})$ with bounded first and second order partial derivatives $\partial_{\hat x}\tilde\Phi,\partial_{\hat x\hat x}\tilde\Phi$.
\end{theorem}

\begin{theorem}\label{theo4}
Assume that {\rm (A1)-(A3)} hold. Then PDE \eqref{pde1} admits a unique solution $\Phi\in C^{1,2}([0,T]\times\mathbb R)$ with bounded $\partial_{\hat x}\Phi$ and $\partial_{\hat x\hat x}\Phi$.
\end{theorem}
\noindent{\bf Proof}\quad
On the one hand, let $\Phi(t,\hat{x})$ be the solution to PDE \eqref{pde1}. Define $\tilde{\Phi}(t,\hat{x})
:=\Phi(t,\tilde\eta_1^{-1}(t)\hat{x})
-\tilde\eta_2(t)\tilde\eta_1^{-1}(t)\hat{x}$. Then we can check that $\tilde{\Phi}(t,\hat{x})$ satisfies PDE \eqref{pde7}. On the other hand, if $\tilde\Phi(t,\hat{x})$ is a solution to PDE \eqref{pde7}, then $\Phi(t,\hat{x}):=
\tilde\Phi(t,\eta_1(t)\hat{x})
+\eta_2(t)\hat {x}$ solves PDE \eqref{pde1}. Therefore, PDE \eqref{pde1} and PDE \eqref{pde7} are equivalent. Using Theorem \ref{thm:tildephi}, the PDE \eqref{pde1} admits a unique solution $\Phi\in C^{1,2}([0,T]\times\mathbb R)$ with bounded $\partial_{\hat x}\Phi$ and $\partial_{\hat x\hat x}\Phi$.
\hfill$\Box$
\begin{theorem}\label{theo5}
Assume that {\rm (A1)-(A3)} hold. Then PDE \eqref{pde2} admits a unique classical solution $U(t,x,\hat{x}):=2P(t)(x-\hat{x})+2\Phi(t,\hat{x})\in C^{1,2,2}([0,T]\times\mathbb R\times\mathbb R)$ with bounded first and second order partial derivatives in $x$ and $\hat x$.
\end{theorem}
\noindent{\bf Proof}\quad
\emph{Existence.} By \eqref{rica1} and \eqref{pde1},
we have:
$$\begin{aligned}
\partial_tU(t,x,\hat{x})
=&2\dot{P}(t)(x-\hat{x})+2\partial_t
\Phi(t,\hat{x})\\
=&-2\big(2AP(t)+Q-R^{-1}B^2P^2(t)
\big)(x-\hat{x})-2\Big\{\partial_{\hat{x}} \Phi(t, \hat{x})\big[A \hat{x}+a\left(\hat{x}
\right)\\
&+(B+\bar{B}) \rho\left(\Phi(t,\hat{x})
\right)\big]+\frac{\sigma_0^2}{2} \partial_{\hat{x}\hat{x}}
\Phi(t, \hat{x})+\frac{1}{2}\dot{q}\left(
\hat{x}\right) +(A+\dot{a}\left(\hat{x}\right)) \Phi(t, \hat{x})+Q\hat{x}\Big\}.
\end{aligned}$$
In addition, we have:
$$\begin{aligned}
&\partial_xU(t,x,\hat{x})=2P(t),\
\partial_{\hat{x}}U(t,x,\hat{x})
=2\partial_{\hat{x}}\Phi(t,\hat{x})-2P(t),\\
&\partial_{xx}U(t,x,\hat{x})=
\partial_{x\hat{x}}U(t,x,\hat{x})
=0,\
\partial_{\hat{x}\hat{x}}
U(t,x,\hat{x})
=2\partial_{\hat{x}\hat{x}}
\Phi(t,\hat{x}).
\end{aligned}
$$
Thus, $U(t,x,\hat{x})
:=2P(t)(x-\hat{x})+2\Phi(t,\hat{x})$ indeed solves PDE \eqref{pde2}.

\emph{Uniqueness.} Since the solution to PDE \eqref{pde2} is the decoupling field of FBSDE \eqref{fbsde1}, then the uniqueness of the solutions to PDE \eqref{pde2} follows from Theorem \ref{theo3}.
\hfill$\Box$
\begin{theorem}
Assume that {\rm (A1)-(A3)} hold. Then the Hamilton-Jacobi equation \eqref{pde3} admits a unique classical solution $V$ with bounded $\partial_{\mu\mu}V$, $\partial_{\mu x}V$ and $\partial_{\mu\mu\mu}V$, and $$\partial_{\mu}V(t,\mu,x)=U(t,x,\int_{\mathbb R}\tilde x \mu(\mathrm{d}\tilde x))$$ for any $(t,\mu,x)\in [0,T]\times\mathcal{P}_2(\mathbb R)\times\mathbb R$ where $U$ is given in Theorem \ref{theo5}.
\end{theorem}
\noindent{\bf Proof}\quad
\emph{Existence.}
Fix $(t_0,\mu)\in [0,T]\times\mathcal{P}_2(\mathbb R)$. We consider the following SDE: for any $\xi_0\in L_{\mathcal{F}_{t_0}^{mf}}^2
      (\Omega^{mf};\mathbb{R})$ such that its law $\mathcal{L}_{\xi_0}=\mu$,
\begin{equation}\label{sde19}
\left\{
\begin{aligned}
\mathrm{d}x^{*,t_0}(t)= & \Big\{A x^{*,t_0}(t)+a\left(
\int_{\mathbb{R}} \tilde{x}\mu_t(\mathrm{d}\tilde{x})\right)
-\frac{1}{2}R^{-1}B^2\big[U(t,
x^{*,t_0}(t),\int_{\mathbb{R}} \tilde{x}\mu_t(\mathrm{d}\tilde{x}))\\
&-\int_{\mathbb{R}}U(t,x,
\int_{\mathbb{R}} \tilde{x}\mu_t(\mathrm{d}\tilde{x}))
\mu_t(\mathrm{d}x)\big]\\
&+(B+\bar{B})\rho\left(\frac{1}{2}
\int_{\mathbb{R}}U(t,
x,\int_{\mathbb{R}} \tilde{x}\mu_t(\mathrm{d}\tilde{x}))
\mu_t(\mathrm{d}x)\right)
\Big\} \mathrm{d}t\\
&+\sigma\mathrm{d}W(t)+\sigma_0 \mathrm{d} W^0(t), \ t \in[t_0, T], \\
x^{*,t_0}(t_0) = &\xi_0,\
\end{aligned}
\right.
\end{equation}
is well-posed on $[t_0,T]$, where $\mu_t:=\mathcal{L}_{x^{*,t_0}(t)
\vert\mathcal{F}^{W^0}_t}$ and $U$ is the unique solution to PDE \eqref{pde2}. Define
\begin{equation}\label{value1}
\begin{aligned}
V(t_0,\mu):=&\mathbb{E}^{mf}\Big\{\int_{t_0}^T\big[Q \int_\mathbb{R}x^2\mu_t(\mathrm{d}x)
+q\left(\int_\mathbb{R}\tilde{x}\mu_t(\mathrm{d}\tilde{x})
\right)
+ R (u^{*,t_0}(t))^2
+r\left( \hat{u}^{*,t_0}(t) \right)\big] \mathrm{d} t  \\
& + G \int_\mathbb{R}x^2\mu_T(\mathrm{d}x)+ g\left(
\int_\mathbb{R}\tilde{x}\mu_T(\mathrm{d}\tilde{x}) \right)\Big\},
\end{aligned}
\end{equation}
where
\begin{equation*}
\begin{aligned}
u^{*,t_0}(t)=&-\frac{1}{2}
R^{-1}B\big[U(t,x,\int_{\mathbb{R}} \tilde{x}\mu_t(\mathrm{d}\tilde{x}))-\int_{\mathbb{R}}U(t,
x,\int_{\mathbb{R}} \tilde{x}\mu_t(\mathrm{d}\tilde{x}))
\mu_t(\mathrm{d}x)\big]\\
&+\rho\left(\frac{1}{2}
\int_{\mathbb{R}}U(t,x,\int_{\mathbb{R}}\tilde{x}\mu_t(\mathrm{d}\tilde{x}))\mu_t(\mathrm{d}x)\right),\\
\hat{u}^{*,t_0}(t)=&\rho\left(
\frac{1}{2}
\int_{\mathbb{R}}U(t,x,\int_{\mathbb{R}}\tilde{x}\mu_t(\mathrm{d}\tilde{x}))
\mu_t(\mathrm{d}x)\right).
\end{aligned}
\end{equation*}
By using formula \eqref{u1} and $U(t,x,\hat{x})=2P(t)(x-\hat{x})+2\Phi(t,\hat{x})$, we have $V(t_0,\mu_0)=\inf_{\alpha\in
\mathcal{U}_{ad}^{mf}[t_0,T]}
J(\mu_0;\alpha(\cdot))$ subject to the state process:
\begin{equation}\label{sde20}
\left\{\begin{aligned}
\mathrm{d}x^{\alpha}(t)= & \left\{A x^\alpha(t)+a\left( \mathbb{E}^{mf}[x^\alpha(t)
\vert\mathcal{F}^{W^0}_t]\right)+B \alpha(t)+\bar{B} \mathbb{E}^{mf}[\alpha(t)
\vert\mathcal{F}^{W^0}_t]
\right\} \mathrm{d}t\\
&+\sigma\mathrm{d} W(t) +\sigma_0 \mathrm{d} W^0(t), \ t \in(t_0, T], \\
x^{\alpha}(t_0)= & \xi_0.
\end{aligned}\right.
\end{equation}
According to the dynamic programming and $U\in C^{1,2,2}([0,T]\times\mathbb R\times\mathbb R)$, we can show that $V\in C^{1,2,2}([0,T]\times\mathcal{P}_2(\mathbb R)\times\mathbb R)$ satisfies the Hamilton-Jacobi equation \eqref{pde3} in the classical sense.

We claim that
\begin{equation}\label{Mouclaim}
\partial_\mu V(t,\mu,x)=U\left(t,x,\int_{\mathbb R}\tilde x\mu(\tilde x)\right).
\end{equation}
We cannot directly differentiate the Hamilton-Jacobi equation \eqref{pde3} in $\mu$ since its solution $V$ is not regular enough. However, the smooth mollifier and a perturbation argument can be applied to overcome this difficult, see \cite[Theorem 3.3]{mou2021wellposedness} and \cite[Theorem 4.2]{cosso2023smooth}. Therefore, without loss of generality, we assume $V$ is smooth. Let us differentiate the Hamilton-Jacobi equation \eqref{pde3} in $\mu$, use \eqref{eq12} and obtain:
\begin{equation}\label{pamuV}
\left\{
\begin{aligned}
&\partial_t(\partial_{\mu}V)(t,\mu,x)+\partial_x
(\partial_{\mu}V)(t,\mu,x)\big[Ax+a\left(\int_{\mathbb R}\tilde x\mu(\mathrm{d}\tilde x)
\right)
-\frac{1}{2}R^{-1}B^2(\partial_{\mu}
V(t,\mu,x)\\
&-\int_{\mathbb R}\partial_{\mu}V(t,\mu,\tilde x)\mu(\mathrm{d}\tilde x))
+(B+\bar{B})\rho\left(\frac{1}{2}\int_{\mathbb R}\partial_{\mu}V(t,\mu,\tilde x)\mu(\mathrm{d}\tilde x)
\right)\big]\\
&+\int_{\mathbb R}\partial_{\mu\mu}V(t,\mu,\tilde x,x)\big[A\tilde x+a\left(\int_{\mathbb R}\bar x\mu(\mathrm{d}\bar x)\right)
-\frac{1}{2}R^{-1}B^2(\partial_{\mu}V(t,\mu,\tilde x)\\
&-\int_{\mathbb R}\partial_{\mu}V(t,\mu,\bar x)\mu(\mathrm{d}\bar x))+(B+\bar{B})\rho\left(\frac{1}{2}\int_{\mathbb R}\partial_{\mu}V(t,\mu,\bar x)\mu(\mathrm{d}\bar x)
\right)\big]\mu(\mathrm{d}\tilde x)\\
&
+\frac{\sigma^2+\sigma_0^2}{2}\partial_{xx\mu}V(t,\mu,x)+\frac{\sigma^2+3\sigma_0^2}{2}\int_{\mathbb R}\partial_{x\mu\mu}V(t,\mu,\tilde x,x)\mu(\mathrm{d}\tilde x)\\
&+\frac{\sigma_0^2}{2}\int_{\mathbb R}\int_{\mathbb R}\partial_{\mu\mu\mu}V(t,\mu,\tilde x,\bar x,x)\mu(\mathrm{d}\tilde x)\mu(\mathrm{d}\bar x)
+\partial_{\mu}V(t,\mu,x)A\\
&+\int_{\mathbb R}\partial_{\mu}V(t,\mu,\tilde x)\mu(\mathrm{d}\tilde x)\dot a\left(\int_{\mathbb R}\tilde x \mu(\mathrm{d}\tilde x)\right)
+2Qx+\dot q\left(\int_{\mathbb R}\tilde x\mu(\mathrm{d}\tilde x)\right)=0,\\
&\partial_{\mu}V(T,\mu,x)=2Gx+\dot{g}
\left(\int_{\mathbb R}x \mu(\mathrm{d}x)\right).
\end{aligned}
\right.
\end{equation}
It can be verified that $\partial_{\mu}V$ is the decoupling field of the MF-FBSDE \eqref{fbsde1}. Since $U(t,x,\hat x)=2P(t)(x-\hat{x})+2\Phi(t,\hat x)$ solves \eqref{pde2}, we can also verify that $\tilde U(t,\mu,x):=U\left(t,x,\int_{\mathbb R}\tilde x\mu(\mathrm{d}\tilde x)\right)$ is the decoupling field of the MF-FBSDE \eqref{fbsde1}. By Theorem \ref{theo3}, we derive the claim \eqref{Mouclaim} holds and $V$ has the desired regularity.

\emph{Uniqueness.}
Let $\bar{V}$ be another classical solution to the Hamilton-Jacobi equation \eqref{pde3}. Again without loss of generality we can assume that $\bar V$ is smooth. Then we can show that $\partial_{\mu}\bar V$ is the decoupling field of the MF-FBSDE \eqref{fbsde1}. By Theorem \ref{theo3}, we can show that $\partial_{\mu}\bar{V}\equiv\partial_{\mu}V$, which implies $V(t,\mu)=\bar V(t,\mu)+c(t)$. Since both $V$ and $\bar V$ are solutions to \eqref{pde3}, we can obtain $\dot c(t)\equiv 0$ and $c(T)=0$. Therefore, $V\equiv \bar V$.
\hfill$\Box$
\section{Solving Problem (LP)}

\subsection{Stochastic Hamiltonian system}

In this subsection, we would like to solve Problem (LP) by the stochastic maximum principle. Firstly, we rewrite Problem (LP) in a vector form. Now $X(\cdot):=(x_1(\cdot),\cdots,x_N(\cdot))^\mathrm{T}$ satisfies the following SDE with the initial condition $\Xi:=(\xi_1,\cdots,\xi_N)^\mathrm{T}$:
\begin{equation}\label{sde4}
\left\{\begin{aligned}
\mathrm{d}X(t)= & \left\{\boldsymbol{A} X(t)+\boldsymbol{a}\left(X(t)\right)
+\boldsymbol{B}\boldsymbol{u}(t)+\boldsymbol{\bar{B}} \left(\boldsymbol{u}(t)\right)\right\} \mathrm{d}t
+\Sigma \mathrm{d} W(t) +\Sigma_0 \mathrm{d} W^0(t), \ t \in(0, T], \\
X(0)= & \Xi,
\end{aligned}\right.
\end{equation}
where $W(t):=(W_1(t),\cdots,W_N(t))^\mathrm{T}$. The corresponding cost functional becomes
\begin{equation}\label{cost2}
\begin{aligned}
\mathcal{J}(\Xi;\boldsymbol{u}(\cdot))= & \mathbb{E}^{lp}\Big\{\int_0^T\left[
X^\mathrm{T}(t)\boldsymbol{Q}X(t)
+\boldsymbol{q}
\left(X(t)\right)
+ \boldsymbol{u}^\mathrm{T}(t)
\boldsymbol{R}\boldsymbol{u}(t)
+\boldsymbol{r}\left( \boldsymbol{u}(t) \right)\right] \mathrm{d} t\\
&+X^\mathrm{T}(T)\boldsymbol{G}X(T)+
 \boldsymbol{g}\left( X(T) \right)\Big\},
\end{aligned}
\end{equation}
where
\begin{equation*}
\begin{aligned}
&\boldsymbol{A}=diag(A,\cdots,A),\
\boldsymbol{B}=diag(B,\cdots,B),\
\Sigma=diag(\sigma,\cdots,\sigma),\
\Sigma_0=(\sigma_0,\cdots,\sigma_0)
^\mathrm{T},\\
&\boldsymbol{Q}=\frac{1}{N}diag(Q,\cdots,Q),\
\boldsymbol{R}=\frac{1}{N}diag(R,\cdots,R),\
\boldsymbol{G}=\frac{1}{N}diag(G,\cdots,G),\\
&\boldsymbol{a}(\cdot),\boldsymbol{b}(\cdot)
:\mathbb{R}^N
\rightarrow\mathbb{R}^N,\\
&\boldsymbol{a}(X)=\left(a(\frac{1}{N}
(1,\cdots,1)X),\cdots,
a(\frac{1}{N}
(1,\cdots,1)X)\right)^\mathrm{T},\\
&\boldsymbol{b}(X)=\left(b(\frac{1}{N}
(1,\cdots,1)X),\cdots,
b(\frac{1}{N}
(1,\cdots,1)X)\right)^\mathrm{T},\\
&\boldsymbol{q}(\cdot),\boldsymbol{r}(\cdot),
\boldsymbol{g}(\cdot):\mathbb{R}^N
\rightarrow\mathbb{R},\\
&\boldsymbol{q}(X)=q(\frac{1}{N}
(1,\cdots,1)X),\
\boldsymbol{r}(X)=r(\frac{1}{N}
(1,\cdots,1)X),\
\boldsymbol{g}(X)=g(\frac{1}{N}
(1,\cdots,1)X).
\end{aligned}
\end{equation*}
Under $\mathrm{(A1)}$, it can be verified that $\boldsymbol{a}(\cdot)
$, $\boldsymbol{q}(\cdot),\boldsymbol{r}(\cdot),$
and $\boldsymbol{g}(\cdot)$ are bounded $C^2$ functions with
bounded 1st and 2nd order derivatives. Hence, we see that under (A1)-(A2), for any $u(\cdot)\in
\mathcal{U}^{lp}_{ad}[0,T]$, SDE \eqref{sde4} admits a unique solution $X(\cdot)\in
L_{\mathbb{F}^{np}}^2
      (\Omega^{np};C(0,T;\mathbb{R}^N))$.

Therefore, Problem (LP) can be reformulated as to minimize cost functional \eqref{cost2} subject to state equation \eqref{sde4} over $\mathcal{U}^{lp}_{ad}[0,T]$. Next, we introduce the adjoint equation as follows
\begin{equation}\label{bsde3}
\left\{\begin{aligned}
\mathrm{d}Y(t)= & -\left[\boldsymbol{A}^\mathrm{T} Y(t)+\dot{\boldsymbol{a}}\left(\bar{X}(t)
\right)^\mathrm{T}Y(t)
+\boldsymbol{Q} \bar{X}(t)
+\frac{1}{2}
\dot{\boldsymbol{q}}
(\bar{X}(t))\right] \mathrm{d}t\\
&+Z(t) \mathrm{d} W(t) +Z_0(t) \mathrm{d} W^0(t), \ t \in[0, T), \\
Y(T)= & \boldsymbol{G}\bar{X}(T)
+\frac{1}{2}
\dot{\boldsymbol{g}}
(\bar{X}(T)).
\end{aligned}\right.
\end{equation}
where $\bar{X}(\cdot)$ is supposed to be the optimal state process of Problem (LP). Assume that $\bar{X}(\cdot)\in L^2_{\mathbb{F}^{lp}}(0,T;\mathbb{R}^d)$. Then BSDE \eqref{bsde3} admits a unique solution $(Y(\cdot),Z(\cdot),Z_0(\cdot))\in
L^2_{\mathbb{F}^{lp}}(\Omega^{lp};C(0,T;
\mathbb{R}^N))\times
L^2_{\mathbb{F}^{lp}}(0,T;\mathbb{R}
^{N\times N})\times
L^2_{\mathbb{F}^{lp}}(0,T;\mathbb{R}
^{N})$. According to the stochastic maximum principle, we have the following optimality condition:
\begin{equation}
\boldsymbol{B}^\mathrm{T}Y(t)
+\boldsymbol{\bar{B}}
^\mathrm{T}Y(t)
+\boldsymbol{R}\bar{u}(t)
+\frac{1}{2}
\dot{\boldsymbol{r}}(\bar{X}(t))
=0.
\end{equation}
Denote
$$Y(\cdot):=(y_1(\cdot),\cdots,
y_N(\cdot))^\mathrm{T},\ Z(\cdot):=\begin{pmatrix}
z_{11}(\cdot) &\cdots &z_{11}(\cdot)   \\
\vdots &\ddots &\vdots    \\
z_{N1}(\cdot) &\cdots &z_{NN}(\cdot)
\end{pmatrix},\
Z_0(\cdot):=(z_{10}(\cdot),\cdots,
z_{N0}(\cdot))^\mathrm{T}.$$
Note that for $X \in \mathbb{R}^{N}$, $u \in \mathbb{R}^N$,
$$\begin{aligned}
\dot{\boldsymbol{a}}(X)=&
\frac{\dot{a}(x^{N})}{N}
(1,\cdots,1)^\mathrm{T}(1,\cdots,1),\
\dot{\boldsymbol{q}}(X)
=\left(\dot{q}(x^{N}),
\cdots,\dot{q}(x^{N})
\right)^\mathrm{T}, \\
\dot{\boldsymbol{r}}(\boldsymbol{u})
=&\left(\dot{r}(u^{N}),
\cdots,\dot{r}(u^{N})
\right)^\mathrm{T}. \\
\end{aligned}$$
Hence, BSDE \eqref{bsde3} has a scalar form. Indeed, for any $1\leq i \leq N$, we get:
\begin{equation}\label{bsde4}
\left\{\begin{aligned}
\mathrm{d}y_i(t)= & -\left[A y_i(t)+\dot{a}
\left(\bar{x}^{N}(t)
\right)y^{N}(t)
+Q \bar{x}_i(t)
+\frac{1}{2}
\dot{q}
\left(\bar{x}^{N}(t)\right)\right] \mathrm{d}t\\
&+\sum_{j=1}^Nz_{ij}(t) \mathrm{d} W_j(t) +z_{i0}(t) \mathrm{d} W^0(t), \ t \in[0, T), \\
y_i(T)= & G\bar{x}_i(T)
+\frac{1}{2}
\dot{g}\left(\bar{x}^{N}(T)\right).
\end{aligned}\right.
\end{equation}
In addition, the optimal control $\bar{u}_i(\cdot)$ of the $i$-th particle satisfies:
\begin{equation}\label{eq6}
By_i(t)+\bar{B}
y^{N}(t)+R\bar{u}_i(t)
+\frac{1}{2}
\dot{r}\left(\bar{u}^{N}(t)\right)=0,\
t\in[0,T].
\end{equation}
Taking the sum from $i$ to $N$ on the above equation and dividing it by $N$, we can obtain
\begin{equation*}
(B+\bar{B})y^{N}(t)
+R\bar{u}^{N}(t)
+\frac{1}{2}
\dot{r}\left(\bar{u}^{N}(t)\right)=0,\ t\in[0,T].
\end{equation*}
Under (A3), it yields that
\begin{equation}\label{eq8}
\begin{aligned}
\bar{u}^{N}(t)=&\rho\left(y^{N}
(t)\right)
=-R^{-1}\big[(B+\bar{B})y^{N}(t)
+\frac{1}{2}
\dot{r}\left(\rho
(y^{N}(t))\right)\big]
\end{aligned}
\end{equation}
where $\rho$ is given in \eqref{eq12}. 
In view of \eqref{eq6} and \eqref{eq8}, the optimal control $\bar{u}_i(\cdot)$ of the $i$-th particle can be solved by
\begin{equation}\label{eq9}
\begin{aligned}
\bar{u}_i(t)=&-R^{-1} B(y_i(t)-y^{N}(t))
+\rho\left(y^{N}(t)\right)\\
=&-R^{-1}\big[By_i(t)
+\frac{1}{2}
\dot{r}\left(\rho
(y^{N}(t))\right)+\bar{B}y^{N}(t)\big],
\ t\in[0,T].
\end{aligned}
\end{equation}
Therefore, system \eqref{sde3}, the corresponding adjoint equation \eqref{bsde4}, along with the average of the states for all particles, can be written as the following stochastic Hamiltonian system at the optimal solution:
\begin{equation}\label{fbsde3}
\left\{\begin{aligned}
\mathrm{d}\bar{x}_i(t)= & \Big\{A \bar{x}_i(t)+a\left
(\bar{x}^{N}(t)\right)
-R^{-1}B^2y_i(t)
-\frac{1}{2}R^{-1}(B+\bar{B})
\dot{r}\left(\rho
(y^{N}(t))\right)\\
&-R^{-1}\bar{B}(2B+\bar{B})y^{N}(t)
\Big\} \mathrm{d}t+\sigma \mathrm{d} W_i(t) +\sigma_0 \mathrm{d} W^0(t), \ t \in(0, T], \\
\mathrm{d}\bar{x}^{N}(t)= & \Big\{A \bar{x}^{N}(t)
+a\left(\bar{x}^{N}(t)\right)
-R^{-1}(B+\bar{B})^2y^{N}(t)\\
&-\frac{1}{2}
R^{-1}(B+\bar{B})\dot{r}\left(\rho
(y^{N}(t))\right)
\Big\} \mathrm{d}t
+\frac{1}{N}\sum_{i=1}^N\sigma \mathrm{d} W_i(t) +\sigma_0 \mathrm{d} W^0(t), \ t \in(0, T], \\
\mathrm{d} y_i(t)= & -\left[A y_i(t)+\dot{a}
\left(\bar{x}^{N}(t)
\right)y^{N}(t)
+Q \bar{x}_i(t)
+\frac{1}{2}
\dot{q}
(x^{N}(t))\right] \mathrm{d}t\\
&+\sum_{j=1}^Nz_{ij}(t) \mathrm{d} W_j(t) +z_{i0}(t) \mathrm{d} W^0(t), \ t \in[0, T), \\
\bar{x}_i(0)=&\xi_i,\ \bar{x}^{N}(0)=  \frac{1}{N}\sum_{i=1}^N\xi_i,\
y_i(T)=  G\bar{x}_i(T)
+\frac{1}{2}
\dot{g}(\bar{x}^{N}(T)).
\end{aligned}\right.
\end{equation}
Now, we would like to decouple the above FBSDE \eqref{fbsde3}.
Similar to the mean field setting, we can derive that
$$y_i(t)=P(t)(\bar{x}_i(t)-\bar{x}^N(t))+\psi(t),\ t\in[0,T],$$
where $P(\cdot)$ is the solution to Riccati equation \eqref{rica1} and $\psi(\cdot)$ solves the following BSDE
\begin{equation}\label{bsde5}
\left\{
\begin{aligned}
\mathrm{d}\psi(t)=&  -\left[A \psi(t)+\dot{a}
\left(\bar{x}^{N}(t)
\right)\psi(t)
+Q \bar{x}^{N}(t)
+\frac{1}{2}
\dot{q}
(\bar{x}^{N}(t))\right] \mathrm{d}t\\
&+\sum_{i=1}^N\Upsilon_{i}(t) \mathrm{d} W_i(t) +\Upsilon_{0}(t) \mathrm{d} W^0(t), \ t \in[0, T), \\
\psi(T)=&G\bar{x}^N(T)+\frac{1}{2}\dot{g}(\bar{x}^{N}(T)).
\end{aligned}
\right.
\end{equation}
Given $\bar{x}^{N}(\cdot)\in L^2_{\mathbb{F}^{lp}}(0,T;\mathbb{R})$, under (A1)-(A3), BSDE \eqref{bsde5} admits a unique strong solution $(\psi(\cdot),\Upsilon_j(\cdot),
1\leq j\leq N,\Upsilon_0(\cdot))
\in L^2_{\mathbb{F}^{lp}}(\Omega^{lp};C(0,T;
\mathbb{R}))\times
L^2_{\mathbb{F}^{lp}}(0,T;\mathbb{R}^N)
\times
L^2_{\mathbb{F}^{lp}}(0,T;\mathbb{R})$.

Furthermore, for $1\leq i \leq N$, $(z_{ij}(\cdot),1\leq j\leq N,z_{i0}(\cdot))$ in BSDE \eqref{bsde4} has the following representation for $1\leq i\leq N$,
$$\begin{aligned}
z_{ii}(t)=&P(t)(\sigma-\frac{\sigma}{N})+\Upsilon_i(t),\
z_{ij}(t)=\Upsilon_j(t)-P(t)\frac{\sigma}{N}, \ \text{for}\  j\neq i,\ 1\leq j\leq N, \\
z_{i0}(t)=&\Upsilon_0(t),
\ t\in[0,T).
\end{aligned}$$
Recalling the second equation in \eqref{fbsde3}, we can derive the following FBSDE
\begin{equation}\label{fbsde4}
\left\{
\begin{aligned}
\mathrm{d}\bar{x}^{N}(t)= & \Big\{A \bar{x}^{N}(t)
+a\left(\bar{x}^{N}(t)\right)
-R^{-1}(B+\bar{B})^2\psi(t)\\
&-
R^{-1}(B+\bar{B})\tilde{r}\left(
\psi(t)\right)
\Big\} \mathrm{d}t
+\frac{1}{N}\sum_{i=1}^N\sigma \mathrm{d} W_i(t) +\sigma_0 \mathrm{d} W^0(t), \ t \in(0, T], \\
\mathrm{d}\psi(t)=&  -\left[A \psi(t)+\dot{a}
\left(\bar{x}^{N}(t)
\right)\psi(t)
+Q \bar{x}^{N}(t)
+\frac{1}{2}
\dot{q}
(\bar{x}^{N}(t))\right] \mathrm{d}t\\
&+\sum_{i=1}^N\Upsilon_{i}(t) \mathrm{d} W_i(t) +\Upsilon_{0}(t) \mathrm{d} W^0(t), \ t \in[0, T), \\
\bar{x}^{N}(0)= & \frac{1}{N}\sum_{i=1}^N\xi_i,\
\psi(T)=G\bar{x}^N(T)+\frac{1}{2}\dot{g}(\bar{x}^{N}(T)).
\end{aligned}
\right.
\end{equation}

The following results are analogues of Theorem \ref{theo2}, Theorem \ref{theo3} and Theorem \ref{theo4}.
\begin{theorem}
Assume that {\rm (A1)-(A3)} hold. Then FBSDE \eqref{fbsde4} admits a unique solution.
\end{theorem}

\begin{theorem}
Assume that {\rm (A1)-(A3)} hold. Then FBSDE \eqref{fbsde3} admits a unique solution.
\end{theorem}
and
\begin{theorem}\label{thm:Nu}
Assume that {\rm (A1)}-{\rm (A3)} hold. Let $P(\cdot)$ and $(\bar{x}^{N}(t),\psi(\cdot),\Upsilon_{i}(\cdot),
1\leq i\leq N,\Upsilon_0(\cdot))$ be the unique solution to \eqref{rica1} and \eqref{fbsde4}, respectively. Then, the optimal control of Problem {\rm (LP)} 
can be written as the following form:
$$\begin{aligned}
\bar{u}_i(t)=&-R^{-1}B
P(t)(\bar{x}_i(t)-\bar{x}^{N}(t))+\rho\left(\psi(t)\right)\\
=&-R^{-1}\big[BP(t)(\bar{x}_i(t)-\bar{x}^{N}(t))+(B+\bar{B})\psi(t)+\tilde{r}\left(\psi(t)\right)\big],\ t\in[0,T],
\end{aligned}
$$
where the corresponding optimal state $\bar{x}_i(\cdot)$ solves
\begin{equation}\label{sde7}
\left\{\begin{aligned}
\mathrm{d}\bar{x}_i(t)= & \Big\{A \bar{x}_i(t)+a\left(\bar{x}^{N}
(t)\right)
-R^{-1}B^2P(t)(\bar{x}_i(t)
-\bar{x}^{N}(t))
+(B+\bar{B})\rho\left(\psi(t)\right)
\Big\} \mathrm{d}t\\
&+\sigma \mathrm{d} W_i(t) +\sigma_0 \mathrm{d} W^0(t), \ t \in(0, T], \\
\bar{x}_i(0)=&\xi_i.
\end{aligned}\right.
\end{equation}
\end{theorem}
\subsection{Hamilton-Jacobi equation}
Analogous to the discussion of Section \ref{sec1}, we investigate the following Hamilton-Jacobi equation corresponding to Problem (LP).
\begin{equation}\label{pde6}
\left\{
\begin{aligned}
&\partial_tV^N(t,X)+\sum_{j=1}^N
\partial_{x_j} V^N(t,X)
\Big[Ax_j+a(x^{N})
+(B+\bar{B})\rho\left(\frac{1}{2}
\sum_{j=1}^N\partial_{x_j}
V^N(t,X)\right)\\
&-\frac{1}{4}R^{-1}B^2(N
\partial_{x_j}V^N(t,X)-
\sum_{j=1}^N
\partial_{x_j} V^N(t,X))
\Big]+\frac{\sigma_0^2}{2}\sum_{k,j=1}^N
\partial_{x_kx_j}V^N(t,X)\\
&+\frac{\sigma^2}{2}\sum_{j=1}^N
\partial_{x_jx_j}V^N(t,X)
+\frac{Q}{N}\sum_{j=1}^Nx_j^2
+q\left(
x^{N}\right)\\
&+R\rho^2\left(
\frac{1}{2N}
\sum_{j=1}^N\partial_{x_j}
V^N(t,X)\right)
+
r\left(\rho\left(
\frac{1}{2N}
\sum_{j=1}^N\partial_{x_j}
V^N(t,X)\right)\right)
=0,\ t\in [0,T)\\
&V^N(T,X)=\frac{G}{N}\sum_{j=1}^Nx_j^2
+g\left(x^{N}\right),
\end{aligned}
\right.
\end{equation}
where $X:=(x_1,\cdots,x_N)^\mathrm{T}$ and $x^{N}:=\frac{1}{N}\sum_{i=1}^Nx_i.$ Similarly, we introduce the following PDE which can be served as the decoupling field of FBSDEs \eqref{fbsde3} and \eqref{fbsde4}, respectively,
\begin{equation}\label{pde5}
\left\{
\begin{aligned}
&\partial_tU_i(t,x_i,\hat x)
+\partial_{x}U_i(t,x_i,\hat x)
\Big[Ax_i+a(\hat x)
-\frac{1}{2}R^{-1}B^2(U_i(t,x_i,
\hat x)-\frac{1}{N}\sum_{j=1}^N
U_j(t,x_j,\hat x))\\
&
+(B+\bar{B})\rho\left(
\frac{1}{2N}
\sum_{j=1}^NU_j(t,x_j,\hat x)\right)
\Big]\\
&+\partial_{\hat x}
U_i(t,x_i,\hat x)\Big[
A\hat x+a(\hat x)+(B+\bar{B})\rho\left(
\frac{1}{2N}
\sum_{j=1}^NU_j(t,x_j,\hat x)\right)\Big]\\
&+\frac{1}{2}\partial_{\hat x\hat x}
U_i(t,x_i,\hat x)
(\frac{\sigma^2}{N}+\sigma^2_0)+
\partial_{x\hat x}
U_i(t,x_i,\hat x)(\frac{\sigma^2}{N}
+\sigma^2_0)+\frac{1}{2}\partial_{xx}
U_i(t,x_i,\hat x)(\sigma^2
+\sigma^2_0)\\
&
+AU_i(t,x_i,\hat x)+\dot{a}(\hat x)\frac{1}{N}
\sum_{j=1}^NU_j(t,x_j,\hat x)
+2Qx_i+\dot{q}\left(\hat x\right)=0,\  t\in[0,T),\\
&U_i(T,x_i,\hat x)
=2Gx_i+\dot{g}\left(\hat x
\right).
\end{aligned}
\right.
\end{equation}
and
\begin{equation}\label{pde4}
\left\{
\begin{aligned}
&\partial_t\Psi(t,\hat x)
+\partial_{\hat x}\Psi(t,\hat x)
\Big[A\hat x+a(\hat x)
+(B+\bar{B})\rho\left(\Psi(t,\hat x)
\right)\Big]\\
&
+\frac{1}{2}
\partial_{\hat x\hat x}
\Psi(t,\hat x)
(\frac{\sigma^2}{N}+\sigma^2_0)
+\frac{1}{2}\dot{q}\left(\hat x\right)
+\big(A+\dot{a}\left(\hat x\right)
\big)
\Psi(t,\hat x)+Q\hat{x}=0,\  t\in[0,T),\\
&\Psi(T,\hat x)=G\hat{x}+\frac{1}{2}\dot{g}\left(
\hat x\right),
\end{aligned}
\right.
\end{equation}
We can follow the same arguments as the ones in Theorems \ref{theo4} and \ref{theo5} to obtain the following theorems for Problem (LP).
\begin{theorem}
Assume that {\rm (A1)-(A3)} hold. Then PDE \eqref{pde4} admits a unique classical solution $\Psi\in C^{1,2}([0,T]\times \mathbb{R})$ with bounded $\partial_{\hat x}\Psi$ and $\partial_{\hat x\hat x}\Psi$.
\end{theorem}

\begin{theorem}
Assume that {\rm (A1)-(A3)} hold. Then PDE \eqref{pde5} admits a unique classical solution $U_i(t,x,\hat x)
:=2P(t)(x-\hat{x})+2\Psi(t,\hat x)\in C^{1,2}([0,T]\times \mathbb{R}\times\mathbb R)$ with bounded first and second order partial derivatives in $x$ and $\hat x$.
\end{theorem}

\begin{theorem}
Assume that {\rm (A1)-(A3)} hold. Then Hamilton-Jacobi equation \eqref{pde6} admits a unique classical solution $V^{N}\in C^{1,2}([0,T]\times \mathbb{R}^N)$ with bounded first and second order derivatives in $X=(x_1,\cdots,x_N)$. Moreover, for $1\leq i\not=j\not=k\leq N$ we have
$$\begin{aligned}
\partial_{x_i}V^N(t,X)
=&\frac{1}{N}U_i(t,x_i,x^{N}),\\
\partial_{x_ix_i}V^N(t,X)
=&\frac{2(N-1)}{N}P(t)+\frac{2}{N^2}
\partial_{\hat x}\Psi(t,x^{N}),\
\partial_{x_ix_j}V^{N}(t,X)
=-\frac{2}{N}P(t)+\frac{2}{N^2}\partial_{\hat x}
\Psi(t,x^{N}),\\
\partial_{x_ix_ix_i}V^{N}(t,X)=&
\partial_{x_ix_ix_j}V^{N}(t,X)=
\partial_{x_ix_jx_k}V^{N}(t,X)=
\frac{2}{N^3}\partial_{\hat x\hat x}
\Psi(t,x^{N}).
\end{aligned}$$
\end{theorem}
\section{Convergence}

In this section, we would like to show that the optimal control problem for the $N$-particle systems converge to the corresponding mean field control problem as the number of particles goes to infinity. In particular, we are able to show the quantitative convergence of the optimal controls of the particle systems as well as their corresponding values. Moreover, we are also able to verify a propagation of chaos property for the associated optimal trajectories and value function.

To begin with, we present several elementary lemmas, which will be useful in sequel.
\begin{lemma}\label{lemma1}
Assume that {\rm (A1)-(A3)} hold. For any $1\leq i\leq N$, we have:
$$\mathbb{E}^{lp}\big[\sup_{0\leq t
\leq T}\vert \bar{x}_i(t)\vert^2\big]\leq K,$$
where $\bar{x}_i(\cdot)$ is the optimal state solving SDE \eqref{sde7}, and $K$ is a constant independent of $N$.
\end{lemma}
\noindent{\bf Proof}\quad
Let $K$ be a constant independent of $N$, which might be different from line to line. By applying the standard BSDE estimation to equation \eqref{bsde5}, we have:
\begin{equation}
  \mathbb{E}^{lp}[\sup_{0\leq t\leq T}\vert\psi(t)\vert^2]\leq K  \mathbb{E}^{lp}[\int_0^T\vert\bar{x}^N(t)\vert^2\mathrm{d}t+1].
\end{equation}
By applying the standard SDE estimation to $\bar{x}^N(\cdot)$ in the FBSDE \eqref{fbsde4}, we have:
\begin{equation*}
\begin{aligned}
\mathbb{E}^{lp}\big[\sup_{0\leq t
\leq T}\vert \bar{x}^{N}(t)\big\vert^2\big]\leq &K\Big\{\mathbb{E}^{lp}\big[\big(
\frac{1}{N}\sum_{i=1}^N\xi_i\big)^2
+\int_0^T(
\vert\bar{x}^{N}(t)\vert^2+\frac{\sigma^2}{N}+\sigma_0^2) \mathrm{d}t\big]+1\Big\}.
\end{aligned}
\end{equation*}
According to Gronwall's inequality and the initial condition, we have:
\begin{equation}\label{eq13}
\begin{aligned}
\mathbb{E}^{lp}\big[\sup_{0\leq t
\leq T}\vert \bar{x}^{N}(t)\big\vert^2\big]\leq K.
\end{aligned}
\end{equation}
We now return to consider SDE \eqref{sde7}. It follows from H\"{o}lder's inequality and the Burkholder-Davis-Gundy (B-D-G) inequality that for any $1\leq i \leq N$,
$$
\begin{aligned}
\mathbb{E}^{lp}\big[\sup_{0\leq t
\leq T}\vert \bar{x}_i(t)\big\vert^2\big]\leq &K\Big\{\mathbb{E}^{lp}\big[\vert\xi_i\vert^2
+\int_0^T\vert
\bar{x}^{N}(t)\vert^2 \mathrm{d}t\big]+1\Big\}\leq K.
\end{aligned}$$
Hence, our conclusion follows.
\hfill$\Box$

Let us consider the following SDE:
\begin{equation}\label{sde9}
\left\{\begin{aligned}
\mathrm{d}x_i^*(t)= & \Big\{A x_i^*(t)+a\left(x^{N,*}
(t)\right)
-R^{-1}B^2P(t)(x_i^*(t)
-\hat{x}(t))
+(B+\bar{B})\rho\left(\varphi(t)\right)
\Big\} \mathrm{d}t\\
&+\sigma \mathrm{d} W_i(t) +\sigma_0 \mathrm{d} W^0(t), \ t \in(0, T], \\
x_i^*(0)=&\xi_i,
\end{aligned}\right.
\end{equation}
where $(\hat{x}(\cdot),\varphi(\cdot),\Lambda_0(\cdot))$ is the strong solution to equations \eqref{fbsde2} and $x^{N,*}(\cdot)=\frac{1}{N}
\sum_{i=1}^Nx_i^*(\cdot)$. We can verify that $x^{N,*}(\cdot)$ solves the following SDE:
\begin{equation}\label{sde10}
\left\{\begin{aligned}
\mathrm{d}x^{N,*}(t)= & \Big\{A x^{N,*}(t)+a\left(x^{N,*}
(t)\right)
-R^{-1}B^2P(t)(x^{N,*}(t)
-\hat{x}(t))+(B+\bar{B})\rho\left(
\varphi(t)\right)
\Big\} \mathrm{d}t\\
&+\frac{1}{N}\sum_{i=1}^N\sigma \mathrm{d} W_i(t) +\sigma_0 \mathrm{d} W^0(t),
 \ t \in(0, T], \\
\bar{x}^{N,*}(0)= & \frac{1}{N}\sum_{i=1}^N\xi_i.
\end{aligned}\right.
\end{equation}
\begin{lemma}\label{lemma2}
Assume that {\rm (A1)-(A3)} hold. Then, for any $1\leq i \leq N$, we obtain:
$$\begin{aligned}
&\mathbb{E}^{lp}\big[\sup_{0\leq t
\leq T}\big\vert x^{N,*}(t)-\hat{x}(t)
\big\vert^2\big]=O\left(\frac{1}{N}
\right),\\
&\mathbb{E}^{lp}\big[\sup_{0\leq t
\leq T}\big\vert \psi(t)-\varphi(t)
\big\vert^2\big]=O\left(\frac{1}{N}
\right),\\
&\mathbb{E}^{lp}\big[\sup_{0\leq t
\leq T}\big\vert \bar{x}^{N}(t)-\hat{x}(t)
\big\vert^2\big]=O\left(\frac{1}{N}
\right),
\end{aligned}
$$
where $x^{N,*}(\cdot)$, $(\bar{x}^{N}(\cdot),\psi(\cdot),\Upsilon_{i}(\cdot),
1\leq i\leq N,\Upsilon_0(\cdot))$, $(\hat{x}(\cdot),\varphi(\cdot),\Lambda_0(\cdot))$ are the strong solutions to \eqref{sde10}, \eqref{fbsde4} and \eqref{fbsde2}, respectively.
\end{lemma}
\noindent{\bf Proof}\quad
By \eqref{sde10} and \eqref{fbsde2}, we have
\begin{equation*}
\left\{\begin{aligned}
\mathrm{d}(x^{N,*}(t)-\hat{x}(t))= & \Big\{A \big[x^{N,*}(t)-\hat{x}(t)\big]
+a\left(x^{N,*}
(t)\right)
-a\left(\hat{x}
(t)\right)\\
&-R^{-1}B^2P(t)\big[x^{N,*}(t)
-\hat{x}(t)\big]
\Big\}\mathrm{d}t
+\frac{1}{N}\sum_{i=1}^N\sigma \mathrm{d} W_i(t),
 \ t \in(0, T], \\
\bar{x}^{N,*}(0)-\hat{x}(0)= & \frac{1}{N}\sum_{i=1}^N\xi_i-x_0.
\end{aligned}\right.
\end{equation*}
It follows from standard SDE estimation that
\begin{equation*}
\begin{aligned}
\mathbb{E}^{lp}\big[\sup_{0\leq t\leq T}\vert x^{N,*}(t)-\hat{x}(t)\vert^2\big]
\leq & K\Big\{\mathbb{E}^{lp}\big[\int_0^T\big(\vert x^{N,*}(t)-\hat{x}(t)\vert^2
+\vert a\left(x^{N,*}
(t)\right)
-a\left(\hat{x}
(t)\right)\vert^2+\frac{\sigma^2}{N}\big)
\mathrm{d}t\\
&+\vert\frac{1}{N}
\sum_{i=1}^N\xi_i-x_0
\vert^2\big]\Big\}.
\end{aligned}
\end{equation*}
Since $\big\{\xi_i\big\}_{1\leq i \leq N}$ is a sequence of i.i.d random variables with finite variance, and $\mathbb{E}[\xi_i]
=x_0$ for any $1\leq i\leq N$, we can derive that
\begin{equation*}
\begin{aligned}
\mathbb{E}^{lp}\big[\vert\frac{1}{N}
\sum_{i=1}^N\xi_i-x_0
\vert^2\big]=&
\frac{1}{N^2}\mathbb{E}^{lp}\Big[
\big\vert
\sum_{i=1}^N(\xi_i-\mathbb{E}^{lp}
[\xi_i])
\big\vert^2\Big]\\
=&\frac{1}{N^2}\Big(\sum_{i=1}^N
\mathbb{E}^{lp}\big[\vert\xi_i-\mathbb{E}^{lp}
[\xi_i]\vert^2\big]+\sum_{j\neq i}^N
\mathbb{E}^{lp}\big[\big(\xi_i-\mathbb{E}^{lp}
[\xi_i]\big)\big(\xi_j-\mathbb{E}^{lp}
[\xi_j]\big)\big]
\Big)\\
=&\frac{1}{N^2}\sum_{i=1}^N
\mathbb{E}^{lp}\big[\vert\xi_i-\mathbb{E}^{lp}
[\xi_i]\vert^2\big]=
O\left(\frac{1}{N}\right).
\end{aligned}
\end{equation*}
Consequently, we can use (A1) to obtain:
\begin{equation*}
\begin{aligned}
\mathbb{E}^{lp}\big[\sup_{0\leq t\leq T}\vert x^{N,*}(t)-\hat{x}(t)\vert^2\big]
\leq & K\mathbb{E}^{lp}\big[\int_0^T\vert x^{N,*}(t)-\hat{x}(t)\vert^2
\mathrm{d}t\big]+O\left(\frac{1}{N}\right).
\end{aligned}
\end{equation*}
In light of Gronwall's inequality, it yields that
$$\mathbb{E}^{lp}\big[\sup_{0\leq t
\leq T}\big\vert x^{N,*}(t)-\hat{x}(t)
\big\vert^2\big]=O\left(\frac{1}{N}
\right).
$$
By the standard BSDE estimates, we have: 
$$\begin{aligned}
&\mathbb{E}^{lp}\left[\sup_{0\leq t
\leq T}\big\vert \psi(t)-\varphi(t)
\big\vert^2+\int_0^T\left(
\sum_{i=1}^N\big\vert\Upsilon_i(t)\big\vert^2+\big\vert \Upsilon_0(t)-\Lambda_0(t)
\big\vert^2\right) \mathrm{d}t
\right]\\
\leq &K
\mathbb{E}^{lp}\big[\vert\bar{x}^{N}(T)
-\hat{x}(T)\vert^2\big].
\end{aligned}$$
Meanwhile, we get:
\begin{equation*}
\left\{\begin{aligned}
\mathrm{d}(\bar{x}^{N}(t)-\hat{x}(t))= & \Big\{A (\bar{x}^{N}(t)
-\hat{x}(t))
+a\left(\bar{x}^{N}(t)\right)
-a\left(\hat{x}(t)\right)\\
&+(B+\bar{B})\big[\rho(\psi(t))-\rho(\varphi(t))\big]\Big\} \mathrm{d}t
+\frac{1}{N}\sum_{i=1}^N\sigma \mathrm{d} W_i(t) ,\ t \in(0, T], \\
\bar{x}^{N}(0)-\hat{x}(0)= & \frac{1}{N}\sum_{i=1}^N\xi_i-x_0,
\end{aligned}\right.
\end{equation*}
which yields that
$$
\begin{aligned}
\mathbb{E}^{lp}\big[\sup_{0\leq t
\leq T}\big\vert \bar{x}^{N}(t)-\hat{x}(t)
\big\vert^2\big]
\leq &K\mathbb{E}^{lp}\big
[\int_0^T\vert \psi(t)-\varphi(t)\vert^2
\mathrm{d}t]+
O\left(\frac{1}{N}
\right)\\
\leq &K\mathbb{E}^{lp}\big
[\int_0^T\big\vert \bar{x}^{N}(t)-\hat{x}(t)
\big\vert^2
\mathrm{d}t]+
O\left(\frac{1}{N}
\right).
\end{aligned}
$$
Using Gronwall's inequality again, we complete the proof.
\hfill$\Box$
\begin{lemma}
Assume that {\rm (A1)-(A3)} hold. Then,
$$\begin{aligned}
&\mathbb{E}^{lp}\big[\sup_{0\leq t
\leq T}\big\vert \Phi(t,\bar{x}^{N}(t))-
\Psi(t,\bar{x}^{N}(t))
\big\vert^2\big]=O\left(\frac{1}{N^2}
\right).
\end{aligned}
$$
where $\Phi$ and $\Psi$ are solutions to \eqref{pde1} and \eqref{pde4}, respectively.
\end{lemma}
\noindent{\bf Proof}\quad
Noted that $\varphi(t)=\Phi(t,\hat{x}(t))$ and $\psi(t)=\Psi(t,\bar{x}^{N}(t))$, $t\in[0,T]$. Applying It\^o's formula to $\Phi(t,\bar{x}^{N}(t))$, we have:
\begin{equation}\label{eq15}
\begin{aligned}
\mathrm{d}\Phi(t,\bar{x}^{N}(t))=&\partial_t
\Phi(t,\bar{x}^{N}(t))\mathrm{d}t+\partial_x
\Phi(t,\bar{x}^{N}(t))
\big[A\bar{x}^{N}(t)+a\left(
\bar{x}^{N}(t)\right)\\
&+(B+\bar{B})\rho\left(\Psi(t,
\bar{x}^{N}(t))\right)\big]\mathrm{d}t+\frac{1}{2}\partial_{xx}
\Phi(t,\bar{x}^{N}(t))
(\frac{\sigma^2}{N}+\sigma^2_0)\mathrm{d}t\\
&+\frac{\sigma}{N}
\sum_{i=1}^N
\partial_x\Phi(t,\bar{x}^{N}(t))\mathrm{d}W_i(t)
+\sigma_0\partial_x\Phi(t,\bar{x}^{N}(t))\mathrm{d}W^0(t),
\end{aligned}
\end{equation}
Setting $\tilde{\varphi}(t):=\Phi(t,\bar{x}^{N}(t))$, $\tilde{\Lambda}_i(t):=\partial_x
\Phi(t,\bar{x}^{N}(t))\frac{\sigma}{N}$, for $1\leq i\leq N$, $\tilde{\Lambda}_0(t):=\partial_x
\Phi(t,\bar{x}^{N}(t))\sigma_0$ for any $t\in[0,T]$. Since $\Phi$ satisfies PDE \eqref{pde1}, combining with \eqref{eq15}, we obtain that $\tilde{\varphi}(\cdot)$ satisfies the following equation:
\begin{equation*}
\left\{
\begin{aligned}
\mathrm{d}\tilde{\varphi}(t)=&-\Big\{
\frac{1}{2}\dot{q}(\bar{x}^{N}(t))+Q\bar{x}^N(t)+\left(A+\dot{a}\left(\bar{x}^{N}(t)
\right)\right)\tilde{\varphi}(t)
-\frac{1}{2}\partial_{xx}\Phi(t,\bar{x}^{N}(t))
\frac{\sigma^2}{N}\\
&+\sigma_0^{-1}\tilde{\Lambda}_0(t)(B+\bar{B})
\big(\rho\left(\tilde{\varphi}(t)\right)-\rho\left(\psi(t)\right)\big)\Big\}\mathrm{d}t\\
&+\sum_{i=1}^N\tilde{\Lambda}_i(t)\mathrm{d}W_i(t)
+\tilde{\Lambda}_0(t)\mathrm{d}W^0(t),\ t\in[0,T),\\
\tilde{\varphi}(T)=
&G\bar{x}^N(T)+\frac{1}{2}\dot{g}(\bar{x}^{N}(T)).
\end{aligned}
\right.
\end{equation*}
According to BSDE \eqref{bsde5} satisfied by $\psi(\cdot)$, we have
\begin{equation*}
\left\{
\begin{aligned}
\mathrm{d}(\tilde{\varphi}(t)-\psi(t))=&-\Big\{
(A+\dot{a}\left(\bar{x}^{N}(t)\right))(\tilde{\varphi}(t)-\psi(t))
-\frac{1}{2}\partial_{xx}\Phi(t,\bar{x}^{N}(t))
\frac{\sigma^2}{N}\\
&+\partial_x\Phi(t,\bar{x}^N(t))(B+\bar{B})
\big(\rho\left(\tilde{\varphi}(t)\right)-\rho\left(\psi(t)\right)\big)\Big\}\mathrm{d}t\\
&+\sum_{i=1}^N(\tilde{\Lambda}_i(t)
-\Gamma_i(t))\mathrm{d}W_i(t)
+(\tilde{\Lambda}_0(t)
-\Gamma_0(t))\mathrm{d}W^0(t),\ t\in[0,T),\\
\tilde{\varphi}(T)-\psi(T)
=&0.
\end{aligned}
\right.
\end{equation*}
By the standard BSDE estimation, it follows that
$$\begin{aligned}
&\mathbb{E}^{lp}\big[\sup_{0\leq t
\leq T}\big\vert \Phi(t,\bar{x}^{N}(t))-
\Psi(t,\bar{x}^{N}(t))
\big\vert^2\big]\leq \frac{K\sigma^4}{N^2},
\end{aligned}
$$
where $K$ depends on the boundedness of $\partial_{x}\Phi$ and $\partial_{xx}\Phi$ in Theorem \ref{theo4}, the time interval $T$, the Lipschitz constant of $\rho$ and the boundedness of $\dot{a}$.
\hfill$\Box$

Define
$$\begin{aligned}
\beta^*(t,x,y):=&-R^{-1}BP(t)(x-y)+ \rho\left(
\Phi(t,y)\right),\\
\bar{\beta}(t,x,y):=&-R^{-1}BP(t)(x-y)+ \rho\left(
\Psi(t,y)\right).
\end{aligned}$$
And then, we have:
$$\bar{u}_i(t)=\bar{\beta}
(t,\bar{x}_i(t),
\bar{x}^{N}(t)),\
u^*(t)=\beta^*(t,x^*(t),\hat{x}(t)).
$$
As a byproduct of the above lemma, we obtain the convergence rate of the closed-loop optimal control for the $N$-particle systems.
\begin{theorem}
Assume that {\rm (A1)-(A3)} hold. Then for any $1\leq i\leq N$,
$$\mathbb{E}^{lp}\big[\sup_{0\leq t
\leq T}\big\vert \beta^*\left(t,\bar{x}_i(t),\bar{x}^{N}(t)\right)
-\bar{\beta}(t,\bar{x}_i(t),\bar{x}^{N}(t))
\big\vert^2\big]=O\left(\frac{1}{N^2}
\right),
$$
where $\bar{x}_i(\cdot)$ solves SDE \eqref{sde7}.
\end{theorem}
The following is the propagation of chaos property for the corresponding optimal trajectories.
\begin{theorem}\label{prop1}
Assume that {\rm (A1)-(A3)} hold. Then
$$\mathbb{E}^{lp}\big[\sup_{0\leq t
\leq T}\big\vert x^*_i(t)-\bar{x}_i(t)
\big\vert^2\big]=O\left(\frac{1}{N}
\right),\ \text{for any $1\leq i\leq N$},
$$
where $\bar{x}_i(\cdot)$ and $x^*_i(\cdot)$ solve \eqref{sde7} and \eqref{sde9}, respectively.
\end{theorem}
\noindent{\bf Proof}\quad
According to the standard SDE estimate and Lemma \ref{lemma2}, we get:
$$\begin{aligned}
&\mathbb{E}^{lp}\big[\sup_{0\leq t
\leq T}\big\vert x^*_i(t)-\bar{x}_i(t)
\big\vert^2\big]\\
\leq &K \mathbb{E}^{lp}\Big\{\int_0^T\big[
\big\vert a\left(
x^{N,*}(t)\right)
-a\left(\bar{x}^{N}(t)
\right)\big\vert^2
+\vert \bar{x}^{N}(t)-\hat{x}(t) \vert^2+\big\vert \rho(\psi(t))-\rho(\varphi(t))\big\vert^2 \big]\mathrm{d}t\Big\}\\
\leq&K\mathbb{E}^{lp}\Big\{\int_0^T\big[
\big\vert x^{N,*}(t)-\hat{x}(t)
\big\vert^2+\big\vert \bar{x}^{N}(t)-\hat{x}(t)
\big\vert^2
+ \big\vert \psi(t)-\varphi(t)
\big\vert^2\big] \mathrm{d}t
\Big\}\\
=&O\left(\frac{1}{N}\right).
\end{aligned}
$$
\hfill$\Box$

Define $\bar{\mu}(X):=\frac{1}{N}\sum_{i=1}^N
\delta_{x_i}$, where $\delta_{x_i}$ denotes the Dirac mass at $x_i\in\mathbb{R}$, $1\leq i \leq N$.
Finally, we show the quantitative convergence of the value functions for the $N$-particle systems. The following theorem presents the asymptotic performance of the cost functional when applied by the optimal control in the mean field control problem.
\begin{theorem}
Assume that {\rm (A1)}-{\rm (A3)} hold. Then,
$$\Big\vert \mathcal{J}(\boldsymbol{u}^*(\cdot))
-\mathcal{J}(\bar{\boldsymbol{u}}(\cdot))\Big\vert
=O\left(\frac{1}{\sqrt{N}}\right)$$
Here,  $\boldsymbol{u}^*:=\left(u^*,\cdots,u^*\right)^\mathrm{T}$ is an 
$N$-vector. $\mathcal{J}(\boldsymbol{u}^*(\cdot))$ represents the cost functional of the $N$-particle systems where each particle applies the optimal control derived from the mean field control problem, and the state process of each particle follows SDE \eqref{sde9}.
\end{theorem}
\noindent{\bf Proof}\quad
Note that
$$\begin{aligned}
\mathcal{J}(\boldsymbol{u}^*(\cdot))=&\frac{1}{N}\sum_{i=1}^N
\mathbb{E}^{lp}\Big\{\int_0^T\big[Qx_i^*(t)^2
+q(x^{N,*}(t))+Ru_i^*(t)^2
+r(u^{N,*}(t))\big] \mathrm{d}t\\
&+Gx_i^*(T)^2
+g(x^{N,*}(T))\Big\}\\
=&\frac{1}{N}\sum_{i=1}^N
\mathbb{E}^{lp}\Big\{\int_0^T\big[Q(x_i^*(t)
-\bar{x}_i(t)+\bar{x}_i(t))^2
+q\left(x^{N,*}(t)\right)\\
&+R(u_i^*(t)
-\bar{u}_i(t)+\bar{u}_i(t))^2
+r\left(u^{N,*}(t)\right) \big]\mathrm{d}t\\
&+G(x_i^*(T)-\bar{x}_i(T)
+\bar{x}_i(T))^2
+g\left(x^{N,*}(T)\right)\Big\}.
\end{aligned}$$
Since $q(\cdot)$, $r(\cdot)$ and $g(\cdot)$ are Lipschitz continuous function, the above equality implies that
$$\begin{aligned}
\mathcal{J}(\boldsymbol{u}^*(\cdot))\leq &\mathcal{J}(\bar{\boldsymbol{u}}(\cdot))
+\frac{1}{N}\sum_{i=1}^N
\mathbb{E}^{lp}\Big\{\int_0^T\big[Q(x_i^*(t)-
\bar{x}_i(t))^2
+2Q\big\vert(x_i^*(t)-
\bar{x}_i(t))\bar{x}_i(t)\big\vert\\
&+K\vert x^{N,*}(t)-\bar{x}^{N}(t)\vert
+R(u_i^*(t)-\bar{u}_i(t))^2
+2R\big\vert (u_i^*(t)-\bar{u}_i(t))\bar{u}_i(t)
\big\vert\\
&+K\vert
u^{N,*}(t)-\bar{u}^{N}(t)\vert\big] \mathrm{d}t
+G(x_i^*(T)
-\bar{x}_i(T))^2\\
&+2G\big\vert(x_i^*(T)
-\bar{x}_i(T))\bar{x}_i(T)\big\vert
+K\vert x^{N,*}(T)-\bar{x}^{N}(T)\vert
\Big\}.
\end{aligned}$$
Further, by the expression of $u^*(\cdot)$ and $\bar{u}(\cdot)$, we obtain:
$$\begin{aligned}
&\Big\vert \mathcal{J}(\boldsymbol{u}^*(\cdot))
-\mathcal{J}(\bar{\boldsymbol{u}}(\cdot))
\Big\vert\\
\leq
&\frac{K}{N}
\sum_{i=1}^N\Big\{\mathbb{E}^{lp}\big[
\sup_{0\leq t\leq T}\vert x_i^*(t)
-\bar{x}_i(t)\vert^2\big]
+\big(\mathbb{E}^{lp}\big[\sup_{0\leq t\leq T}\vert \bar{x}_i(t)\vert^2\big]\big)
^{\frac{1}{2}}\big(\mathbb{E}^{lp}
\big[\sup_{0\leq t\leq T}\vert
x_i^*(t)-\bar{x}_i(t)\vert^2\big]\big)
^{\frac{1}{2}}\\
&+\mathbb{E}^{lp}\big[
\sup_{0\leq t\leq T}\vert x^{N,*}(t)
-\bar{x}^{N}(t)\vert\big]
+\mathbb{E}^{lp}\big[
\sup_{0\leq t\leq T}\vert x^{N,*}(t)
-\hat{x}(t)\vert\big]\\
&+\mathbb{E}^{lp}\big[
\sup_{0\leq t\leq T}\vert \bar{x}^{N}(t)
-\hat{x}(t)\vert\big]+\mathbb{E}^{lp}\big[
\sup_{0\leq t\leq T}\vert \psi(t)
-\varphi(t)\vert\big]\\
&
+\big(1+\mathbb{E}^{lp}\big[\sup_{0\leq t\leq T}\vert\bar{x}_i(t)\vert^2
\big]+\mathbb{E}^{lp}\big[\sup_{0\leq t
\leq T}\vert\bar{x}^{N}(t)\vert^2
\big]\big)^{\frac{1}{2}}
\big(\mathbb{E}^{lp}\big[\sup_{0\leq t
\leq T}\vert x_i^*(t)-\bar{x}_i(t)\vert^2\\
&+\sup_{0\leq t
\leq T}\vert\bar{x}^{N}(t)
-\hat{x}(t)\vert^2
+\sup_{0\leq t
\leq T}\vert\psi(t)
-\varphi(t)\vert^2\big]\big)
^{\frac{1}{2}}\Big\}.
\end{aligned}
$$
Finally, based on Lemma \ref{lemma1}, \ref{lemma2} and Theorem \ref{prop1}, we can draw the conclusion.
\hfill$\Box$
\begin{theorem}
Assume that {\rm (A1)-(A3)} hold. Then
$$\Big\vert V(t,\bar{\mu}(\bar{X}(t)))
-V^N(t,\bar{X}(t))\Big\vert
=O\left(\frac{1}{N}\right),$$
where $\bar{X}(\cdot):=(\bar{x}_1(\cdot),
\cdots,\bar{x}_N(\cdot))^\mathrm{T}$, $\bar{x}_i(\cdot)$ solves SDE \eqref{sde7}, for $1\leq i \leq N$.
\end{theorem}
\noindent{\bf Proof}\quad
According to the definition of value function, we have:
$$
\begin{aligned}
&\big\vert V(t,\bar{\mu}(\bar{X}(t)))
-V^N(t,\bar{X}(t))
\big\vert\\
=&\mathbb{E}^{lp}\Big\{\int_t^T\big[R\rho^2\left(\Phi(s,\bar{x}^{N}(s))
\right)-R\rho^2\left(\Psi(s,\bar{x}^{N}(s))
\right)+r\left(\rho\left(
\Phi(s,\bar{x}^{N}
(s))\right)\right)-r\left(\rho\left(
\Psi(s,\bar{x}^{N}
(s))\right)\right)\big]\mathrm{d}s\Big\}\\
\leq&\Big(\mathbb{E}^{lp}\Big\{\int_t^TR\big\vert\rho\left(\Phi(s,\bar{x}^{N}(s))
\right)+\rho\left(\Psi(s,\bar{x}^{N}(s))
\right)
\big\vert^2\mathrm{d}s\big\}\Big)^{\frac{1}{2}}\cdot\\
&\Big(\mathbb{E}^{lp}\Big\{\int_t^TRK^2\left\vert
\Phi(s,\bar{x}^{N}(s))
-\Psi(s,\bar{x}^{N}(s))
\right\vert^2\mathrm{d}s\Big\}\Big)^{\frac{1}{2}}+K\mathbb{E}^{lp}\Big\{\int_t^T\left\vert
\Phi(s,\bar{x}^{N}(s))
-\Psi(s,\bar{x}^{N}(s))
\right\vert \mathrm{d}s\Big\}\\
\leq & K\Big(\mathbb{E}^{lp}\Big\{\int_t^T\big[\rho^2(\Phi(s,0))+\rho^2(\Psi(s,0))
+\vert\bar{x}^{N}(s)\vert^2\big]\mathrm{d}s\Big\}\Big)^{\frac{1}{2}}\cdot\\
&\Big(\mathbb{E}^{lp}\big\{\int_t^T\left\vert
\Phi(s,\bar{x}^{N}(s))
-\Psi(s,\bar{x}^{N}(s))
\right\vert^2\mathrm{d}s\big\}\Big)^{\frac{1}{2}}+K\mathbb{E}^{lp}\Big\{\int_t^T\left\vert
\Phi(s,\bar{x}^{N}(s))
-\Psi(s,\bar{x}^{N}(s))
\right\vert \mathrm{d}s\Big\}\\
\leq &K \Big(\mathbb{E}^{lp}\Big\{\int_t^T\big\vert
\Phi(s,\bar{x}^{N}(s))
-\Psi(s,\bar{x}^{N}(s))\big\vert^2 \mathrm{d}s\Big\}\Big)^{\frac{1}{2}}\\
=& O(\frac{1}{N}),
\end{aligned}
$$
where the last inequality is obtained based on \eqref{eq13}, and the continuity of
$\rho$, $\Phi$ and $\Psi$.
\hfill$\Box$
\section*{Appendix}
\setcounter{equation}{0}
\setcounter{subsection}{0}
\renewcommand{\theequation}{A.\arabic{equation}}
\renewcommand{\thesubsection}{A.\arabic{subsection}}
\renewcommand\thetheorem{A.\arabic{theorem}}
\renewcommand\thelemma{A.\arabic{lemma}}
\noindent{\bf Proof of Theorem \ref{thm:tildephi}.}\quad
From Theorem 2.6 in \cite{delarue2002existence}, FBSDE \eqref{fbsde8} admits a unique solution $(\hat{x}^{\eta,t_0,x_0}(\cdot),
\varphi^{\eta,t_0,x_0}(\cdot)$, $
\Lambda_0^{\eta,t_0,x_0}(\cdot))$ for any $(t_0,x_0)\in[0,T]\times\mathbb{R}$.
Next, we establish the well-posedness of the PDE \eqref{pde7}. To do so, we consider the following FBSDE which can be viewed as a formal differentiation of FBSDE \eqref{fbsde8} with respect to $x_0$.
\begin{equation}\label{fbsde6}
\left\{
\begin{aligned}
\mathrm{d} \nabla\hat{x}^{\eta,t_0,x_0}(t)= & \Big\{-R^{-1}(B+\bar{B})^2\big[\tilde{\eta}_1^{t_0,x_0}(t)\nabla\varphi^{\eta,t_0,x_0}(t)+\nabla\tilde{\eta}_1^{t_0,x_0}(t)\varphi^{\eta,t_0,x_0}(t)\big]\\
&+
\big[a^\eta\left(t,
\hat{x}^{\eta,t_0,x_0}(t)\right)
-R^{-1}(B+\bar{B})\tilde{r}^\eta\left(t,
\hat{x}^{\eta,t_0,x_0}(t),
\varphi^{\eta,t_0,x_0}(t)\right)
\big]\nabla\tilde{\eta}_1^{t_0,x_0}(t)\\
&+
\big[\partial_x a^\eta\left(t,
\hat{x}^{\eta,t_0,x_0}(t)\right)
\nabla\hat{x}^{\eta,t_0,x_0}(t)\\
&-R^{-1}(B+\bar{B})\big(\partial_x\tilde{r}
\left(t,\hat{x}^{\eta,t_0,x_0}(t),
\varphi^{\eta,t_0,x_0}(t)\right)
\nabla\hat{x}^{\eta,t_0,x_0}(t)\\
&+\partial_\varphi\tilde{r}
\left(t,\hat{x}^{\eta,t_0,x_0}(t),
\varphi^{\eta,t_0,x_0}(t)
\right)\nabla\varphi^{\eta,t_0,x_0}(t)\big)\big]\tilde{\eta}_1^{t_0,x_0}(t)
\Big\} \mathrm{d}t\\
&+\nabla\tilde{\eta}_1^{t_0,x_0}(t)\sigma_0\mathrm{d}W^0(t), \ t \in[t_0, T], \\
\mathrm{d}\nabla\varphi^{\eta,t_0,x_0}(t)
=&-\Big\{\big[
A+\tilde{a}^\eta\left(t,
\hat{x}^{\eta,t_0,x_0}(t)
\right)-R^{-1}(B+\bar{B})^2\tilde{\eta}_2^{t_0,x_0}(t)\big]
\nabla\varphi^{\eta,t_0,x_0}(t)\\
&+\big[
\partial_x\tilde{a}^\eta\left(t,
\hat{x}^{\eta,t_0,x_0}(t)
\right)
\nabla\hat{x}^{\eta,t_0,x_0}(t)-R^{-1}(B+\bar{B})^2\nabla\tilde{\eta}_2^{t_0,x_0}(t)
\big]
\varphi^{\eta,t_0,x_0}(t)\\
&+\nabla\tilde{\eta}_2^{t_0,x_0}(t)\big[a^\eta\left(t,
\hat{x}^{\eta,t_0,x_0}(t)\right)
-R^{-1}(B+\bar{B})\tilde{r}^\eta
\left(t,\hat{x}^{\eta,t_0,x_0}(t),
\varphi^{\eta,t_0,x_0}(t)\right)\big]\\
&
+\tilde{\eta}_2^{t_0,x_0}(t)
\big[\partial_xa^\eta\left(t,
\hat{x}^{\eta,t_0,x_0}(t)\right)
\nabla\hat{x}^{\eta,t_0,x_0}(t)\\
&-R^{-1}(B+\bar{B})\big(\partial_x\tilde{r}
\left(t,\hat{x}^{\eta,t_0,x_0}(t),
\varphi^{\eta,t_0,x_0}(t)\right)
\nabla\hat{x}^{\eta,t_0,x_0}(t)\\
&+\partial_\varphi\tilde{r}
\left(t,\hat{x}^{\eta,t_0,x_0}(t),
\varphi^{\eta,\eta,t_0,x_0}(t)
\right)\nabla\varphi^{\eta,t_0,x_0}(t)\big)\\
&+\partial_x\tilde{q}^\eta\left(t,
\hat{x}^{\eta,t_0,x_0}(t)\right)
\nabla\hat{x}^{\eta,t_0,x_0}(t)
\Big\}\mathrm{d}t+\nabla\Lambda_0^{\eta,t_0,x_0}(t)
\mathrm{d}W^0(t), \ t\in[t_0,T],\\
\nabla\hat{x}^{\eta,t_0,x_0}(t_0)= &1,\
\nabla\varphi^{\eta,t_0,x_0}(T)
=\partial_x\tilde{g}^\eta
(T,\hat{x}^{\eta,t_0,x_0}(T))
\nabla\hat{x}^{\eta,t_0,x_0}(T).
\end{aligned}
\right.
\end{equation}
Here, $\nabla\tilde{\eta}_1^{t_0,x_0}
(\cdot)$ and $\nabla\tilde{\eta}_2^{t_0,x_0}
(\cdot)$ satisfy the following linear equations
\begin{equation}\label{sde17}
\left\{
\begin{aligned}
\mathrm{d} \nabla\tilde{\eta}_1^{t_0,x_0}(t)= 
&-\Big\{\big[A-R^{-1}(B+\bar{B})^2
\tilde{\eta}_2^{t_0,x_0}(t)\big]\nabla\tilde{\eta}_1^{t_0,x_0}(t)\\
&-R^{-1}(B+\bar{B})^2
\nabla\tilde{\eta}_2^{t_0,x_0}(t)\tilde{\eta}_1^{t_0,x_0}(t)\Big\}
\mathrm{d}t,\ t\in[t_0,T],\\
\nabla\tilde{\eta}_1^{t_0,x_0}(T)= &0,
\end{aligned}
\right.
\end{equation}
and
\begin{equation}\label{sde18}
\left\{
\begin{aligned}
\mathrm{d} \nabla\tilde{\eta}_2^{t_0,x_0}(t)= & -\Big\{
\big[2A+\tilde{a}^\eta\left(t,\hat{x}^{\eta,t_0,x_0}(t)\right)-2R^{-1}(B+\bar{B})^2\tilde{\eta}_2^{t_0,x_0}(t)\big]
\nabla\tilde{\eta}_2^{t_0,x_0}(t)\\
&+\partial_x\tilde{a}^\eta\left(t,
\hat{x}^{\eta,t_0,x_0}(t)\right)
\nabla\hat{x}^{\eta,t_0,x_0}(t)
\tilde{\eta}_2^{t_0,x_0}(t)
\Big\}\mathrm{d}t,\ t\in[t_0,T],\\
\nabla\tilde{\eta}_2^{t_0,x_0}(T)= &0,
\end{aligned}
\right.
\end{equation}
which are well-posed on $[t_0,T]$.

\emph{Existence.} Let us consider
FBSDE \eqref{fbsde6} for any $(t_0,x_0)\in[0,T]\times\mathbb{R}$.
Under (A1), we can assume that the coefficients in FBSDE \eqref{fbsde6} are bounded by $
\tilde{M}$.
According to Corollary 2.8 in \cite{delarue2002existence}, we have $\vert\partial_{\hat x}\tilde{\Phi}\vert\leq
\tilde{M}$ (if needed we enlarge $\tilde M$). From the standard theory of FBSDEs (see \cite{zhang2017backward}), there exists a $\delta>0$ depending on $\tilde{M}$ such that FBSDE \eqref{fbsde6} admits a unique solution on $[T-\delta,T]$. Using $\vert\partial_{\hat x}\tilde{\Phi}\vert\leq
\tilde{M}$, FBSDE \eqref{fbsde6} is well-posed on $[T-2\delta,T-\delta]$ with the terminal condition changed to $\nabla\varphi^{\eta,t_0,x_0}(T-\delta)
=\partial_{\hat x}\tilde{\Phi}
(T-\delta,\hat{x}^{\eta,t_0,x_0}(T-\delta))
\nabla\hat{x}^{\eta,t_0,x_0}(T-\delta)$. Finally, repeating finite times, we should be able to obtain that
FBSDE \eqref{fbsde6} admits a unique strong solution $(\nabla\hat{x}^{\eta,t_0,x_0}(\cdot),
\nabla\varphi^{\eta,t_0,x_0}(\cdot),
\nabla\Lambda_0^{\eta,t_0,x_0}(\cdot))\in L_{\mathbb{F}^0}^2
      (\Omega;C(0,T;\mathbb{R}))\times L_{\mathbb{F}^0}^2
      (\Omega;C(0,T;\mathbb{R}))\times L^2_{\mathbb{F}^0}(0,T
      ;\mathbb{R})
$ for any $(t_0,x_0)\in[0,T]\times\mathbb{R}$ and $\partial_{\hat x}\tilde{\Phi}\in C_b([0,T]\times\mathbb{R})$. By further differentiating \eqref{fbsde6} in $x_0$ and using $\vert\partial_{\hat x}\tilde{\Phi}\vert\leq
\tilde{M}$, we can show that $\partial_{\hat x\hat x}\tilde{\Phi}\in C_b([0,T]\times\mathbb{R})$. With the continuity of $\partial_{\hat x}\tilde{\Phi}$ and $\partial_{\hat x\hat x}\tilde{\Phi}$, it is standard to verify $\tilde{\Phi}\in C^{1,2}([0,T]\times \mathbb{R})$ and it satisfies PDE \eqref{pde7}.

\emph{Uniqueness.} Let $\bar{\tilde{\Phi}}\in C^{1,2}([0,T]\times\mathbb{R})$ be another classical solution to PDE \eqref{pde7} with bounded $\partial_{x}\bar{\tilde{\Phi}}$ and $\partial_{xx}\bar{\tilde{\Phi}}$. Introduce the following SDE
\begin{equation}\label{sde14}
\left\{
\begin{aligned}
\mathrm{d}\bar{\hat{x}}^{\eta,t_0,x_0}(t)= & \Big\{-R^{-1}(B+\bar{B})^2\tilde{\eta}_1^{t_0,x_0}(t)
\bar{\tilde{\Phi}}(t,
\bar{\hat{x}}^{\eta,t_0,x_0}(t))+\big[a^\eta\left(t,
\bar{\hat{x}}^{\eta,t_0,x_0}(t)\right)\\
&-R^{-1}(B+\bar{B})\tilde{r}^\eta\left(t,
\bar{\hat{x}}^{\eta,t_0,x_0}(t),
\bar{\tilde{\Phi}}(t,
\bar{\hat{x}}^{\eta,t_0,x_0}(t))\right)\big]\tilde{\eta}_1^{t_0,x_0}(t)\Big\} \mathrm{d}t\\
&+\tilde{\eta}_1^{t_0,x_0}(t)\sigma_0\mathrm{d} W^0(t), \ t \in[t_0, T],  \\
\bar{\hat{x}}^{\eta,t_0,x_0}(t_0)= & x_0,
\end{aligned}
\right.
\end{equation}
Here, $\tilde{\eta}_1^{t_0,x_0}(\cdot)$ and
$\tilde{\eta}_2^{t_0,x_0}(\cdot)$
satisfy equations \eqref{sde15} and \eqref{sde16} respectively with $(\hat{x}(\cdot),\varphi(\cdot))$
replaced by
$\left((\tilde{\eta}_1^{t_0,x_0})^{-1}(\cdot)
\bar{\hat{x}}^{\eta,t_0,x_0}(\cdot)\right.$,
$\left.\bar{\tilde{\Phi}}(\cdot,
\bar{\hat{x}}^{\eta,t_0,x_0}(\cdot))
+\tilde{\eta}_2^{t_0,x_0}(\cdot)
(\tilde{\eta}_1^{t_0,x_0})^{-1}
(\cdot)
\bar{\hat{x}}^{\eta,t_0,x_0}(\cdot)
\right)$. It is clear that \eqref{sde14} is well-posed on $[t_0,T]$.
Define for any $t\in[t_0,T]$, $$\bar{\varphi}^{\eta,t_0,x_0}(t)
:=\bar{\tilde\Phi}(t,
\bar{\hat{x}}^{\eta,t_0,x_0}(t)),\
\bar{\Lambda}^{\eta,t_0,x_0}_0(t):=
\partial_x\bar{\tilde\Phi}(t,
\bar{\hat{x}}^{\eta,t_0,x_0}(t))
\tilde{\eta}_1^{t_0,x_0}(t)\sigma_0.$$
Noted that $\bar{\tilde\Phi}$ is a solution to PDE \eqref{pde7}, we have $(\bar{\hat{x}}^{\eta,t_0,x_0}(\cdot),
\bar{\varphi}^{\eta,t_0,x_0}(\cdot),
\bar{\Lambda}^{\eta,t_0,x_0}_0(\cdot))$ solves FBSDE \eqref{fbsde8}. In view of the uniqueness of FBSDE \eqref{fbsde8}, it follows that $\tilde{\Phi}
=\bar{\tilde\Phi}$.
\hfill$\Box$
\section*{Acknowledgements}
The first author was supported by the China Scholarship Council (No. 202206220077).
The second author is supported by CityU Start-up Grant 7200684, Hong Kong RGC Grant ECS
21302521, Hong Kong RGC Grant GRF 11311422 and Hong Kong RGC Grant GRF 11303223.
The third author is supported by the National Key Research and Development Program of China
under Grant 2023YFA1009200, the National Natural Science Foundation of China (No. 11831010,
61961160732), the Natural Science Foundation of Shandong Province (No. ZR2019ZD42), the Taishan
Scholars Climbing Program of Shandong (No. TSPD20210302). The fourth author is supported
by Singapore MOE (Ministry of Education) AcRF Grant A-8000453-00-00, IoTex Foundation Industry
Grant A-8001180-00-00 and NSFC Grant No. 11871364.

\end{document}